# The Ising model in our dimension and our times

**Andrei Okounkov**

**Abstract**

While the author is a professional mathematician, he is by no means an expert in the subject area of these notes. The goal of these notes is to share the author's personal excitement about some results of Hugo Duminil-Copin with mathematics enthusiasts of all ages, using maximally accessible, yet precise mathematical language. No attempt has been made to present an overview of the current state field, its history, or to place this narrative in any kind of broader scientific or social context. See the references in Section 5 for both professional surveys and popular science accounts that will certainly give the reader a broader and deeper understanding of the material.

Phase transitions are dramatic physical phenomena. A physical system undergoing a phase transition may exhibit large spatial fluctuations, a detailed understanding of which presents an important challenge to physicists and mathematicians alike. Thanks to Hugo Duminil-Copin and his collaborators, the recent years saw a great progress in our understanding of the phase transition in the 3-dimensional Ising model, perhaps the most famous model of mathematical statistical physics.

The goal of these notes is to explain an introductory portion of this progress to the broadest possible audience of mathematics enthusiasts. Before we get to say anything of substance about the new results, there is a certain amount of language to develop and background to review.

## 1. Mathematics and physics

Mathematics provides the universal language of science. While human languages have words that describe natural phenomena, they lag far, far behind the language of mathematics in their precision and predictive power. It is easy to fill a sizable volume with quotes to this effect from the most prominent scientists of all epochs[1].

Wouldn't the task of writing these notes be really simple if mathematics were *only* a language ? There would probably be usable automatic translation available at a click. In fact, it is a very common request to translate from mathematics to a natural language[2]. Richard Feynman, in particular, talks about it in the second of his 1964 lectures about the Character of Physical Law[3].

While mathematics has its special words and symbols, as well as grammatical rules that govern what is the correct logical use of these symbols and what is not, the real treasure of mathematics is the much deeper level of understanding that this language empowers. By defining the boundaries of precise reasoning, and removing all other boundaries between ideas, mathematics allows humans to use the most inventive and unexpected mathematical constructions and arguments to discover deep truths about the world around us. Instead of being lost in the confusing woods of natural languages, mathematics makes it possible for our thought to fly *safely*.

These notes are about mathematical physics, the field where mathematics and physics come together. A mathematical physicist starts by *defining* her or his object of study, introducing a mathematical object that captures some essential features of one or many physical phenomena. One calls it a *model*, which, unlike many other words used by mathematicians, is a term that stays fairly close to its meaning in natural languages. Having defined a model, a mathematical physicist is free to be arbitrarily creative in her or his choice of mathematical tools to study this model. This investigation is going through all the natural stages of research in mathematics: one asks precise questions, considers examples, formulates conjectures, obtains partial results, and, as a proof of having achieved a really good understanding of the model, one can prove mathematical theorems about it.

---

1 Gibbs measures, named so in honor of J. Willard Gibbs (1839–1903), will play the central role in our narrative. Gibbs is remembered as very unsociable and the only words he ever said in the Yale faculty meeting were *Mathematics is a language*. See the biography [30] of Gibbs by Muriel Rukeyser. She also wrote a poem about Gibbs inspired by this quote.

2 In the narrator's personal experience, good progress in science often happens when trying to answer the opposite question, namely, *can you translate what you just said to mathematics* ?

3 Feynman says, in particular, this about translating mathematics: *But I do not think it is possible, because mathematics is not just another language. Mathematics is a language plus reasoning; it is like a language plus logic. Mathematics is a tool for reasoning. It is in fact a big collection of the results of some person's careful thought and reasoning. By mathematics it is possible to connect one statement to another.*

For example, in his Mathematical Principles of Natural Philosophy, Newton introduced differential equation as a mathematical language to describe the motion of celestial as well as terrestrial bodies. This gives him a *model* for motion of planets around the sun. In the approximation that ignores the mutual attraction of the planets, he then mathematically *proves* the planets follow Kepler's empirical laws of planetary motion.

The language and the models evolve. Each chapter in that great book of the Universe to which Galileo refers[4] in Il Saggiatore is written in a new mathematical language that has to be discovered every time. Newtonian mechanics is an approximation that is good at modeling some phenomena but not others. Quantum mechanics had to be created to describe the behavior of molecules, atoms, and other tiny constituents of the universe. Statistical physics had to be created to describe phenomena in which the myriads of particles that form planets and other macroscopic objects don't just move as one, but instead create very complex patterns and materials through spatial interactions. The actual mathematics used in each case is very different. The Ising model, which will be our focus of attention in this narrative, is perhaps the most famous model of statistical physics.

A question often asked about mathematical physics is: where is the boundary between mathematics and physics in it ? In the personal view of this narrator, there is no boundary[5]. It is a really joint endeavor between mathematics and physics, where each side contributes something extremely important. Among other things, physics provides invaluable *intuition*, rooted in laboratory and numerical experiments, as well as parallels and correspondences that extends across different branches of physics. These can guide mathematics at any of the research stages discussed above. For mathematical physicists, following the logic of the subject is much more important than departmental affiliation. For instance, the first truly amazing mathematical result about the Ising model was obtained by Lars Onsager, the winner of the 1968 Nobel Prize in chemistry. In the narrator's personal experience, physicists are very proud when they find a mathematical proof and mathematicians are very proud when they discover a good physical explanation.

---

**4** *Philosophy is written in this grand book, which stands continually open before our eyes (I say the 'Universe'), but can not be understood without first learning to comprehend the language and know the characters as it is written. It is written in mathematical language, and its characters are triangles, circles and other geometric figures, without which it is impossible to humanly understand a word; without these one is wandering in a dark labyrinth.*

**5** If it really exists, the boundary is as diffuse as the boundary in Figure (1)

It is fitting to end this quotation-filled section with a quote from the hero of these notes Hugo Duminil-Copin: "*Mathematical Physics gathers everything I always dreamt of as a researcher*: *it satisfies my curiosity to understand the physical world in which we live, and it rewards the mathematicians with beautiful and elegant rigorous proofs*."

## 2. The Ising model

### 2.1. Stuff fluctuates in space

It is good to have a picture in mind as we discuss the definition of the Ising model. Here is one, a simulation by Stanislav Smirnov. Its meaning will be made clear gradually.

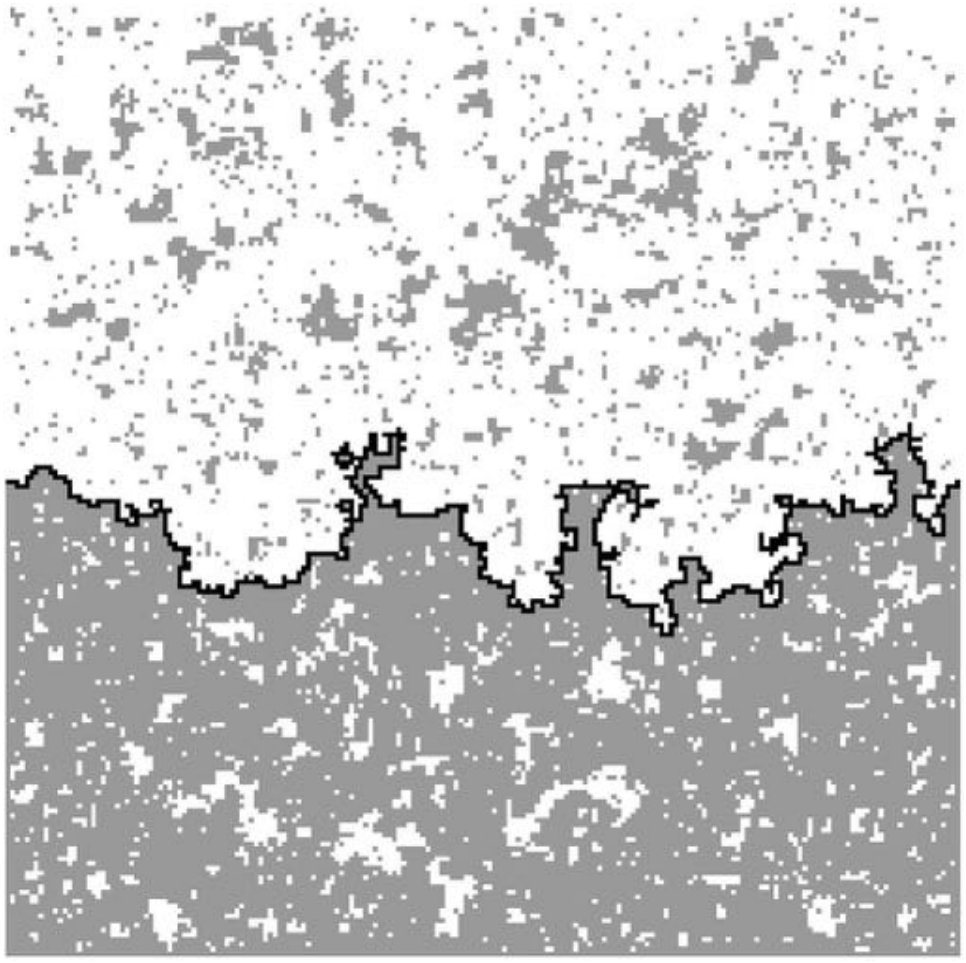

(1)

Clearly, this is something *random*. The language of statistical physics is based on randomness and probabilities.

Our world is fundamentally random. In statistical physics randomness is introduced from the very beginning[6]. From the very beginning, statistical physics talks about the *probabilities* for a physical system to be in such or such state.

Very importantly, the randomness in Figure (1) happens in *space*, here a 2-dimensional space. In other words, in (1) we have a random spatial pattern. Note that while obviously complex, this pattern is not pure noise. We see a very diffuse boundary between black and white, with many islands or lakes of one color inside another. These have intricate shapes and may be nested, that is, there can be an island on a lake in the middle of a larger island on a larger lake, et cetera.

People are usually introduced to probability theory through coin tosses, rolls of dice, and similar random events that have a few possible outcomes and no spatial structure.

---

6 In quantum mechanics, randomness is also present from the very beginning. In principle, Newtonian mechanics makes exact predictions about the behaviour of its models. However, for systems of large size and complexity, think Avogadro number many billiard balls bouncing off each other, these predictions are so complicated as to be effectively random. This is the subject of ergodic theory, the development of which was very much stimulated by the quest to see statistical physics emerge from Newtonian mechanics.

Successive games of chance and similar data sets (think stock prices, air temperatures etc.) produce random time series like the one in Figure (2). These have a 1-dimensional structure to it. They are like beads threaded by the axis of time. In probability theory, these are known as random processes.

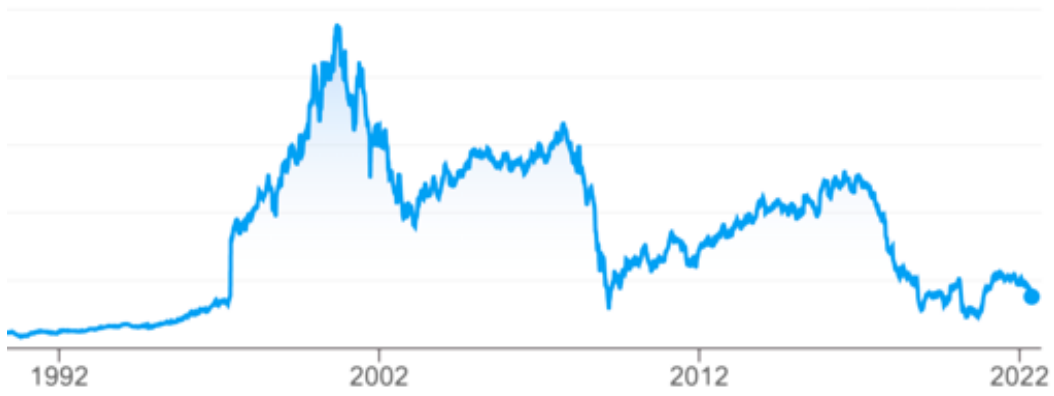

(2)

Statistical physics really starts in dimensions two or more, studying random objects fluctuating in the corresponding number of dimensions. Importantly, the behavior of the Ising model (and most models of statistical physics) very strongly depends on the dimension. The Ising model is dull in dimension 1, very interesting in dimensions 2 and 3, and regresses to the generic, and hence not as exciting[7], Gaussian behaviour in dimensions $\geq 4$.

Many past glorious mathematical successes of statistical physics concern the 2-dimensional Ising model. While we will say a few words about it, our goal in this narrative is to report on the great recent progress in our physical dimension 3 achieved by Hugo Duminil-Copin and his collaborators. This will come in due course. We have not defined the Ising model yet.

### 2.2. A lattice in space

While it is good to imagine that a random process like the one in (2) happens in continuous time, the actual data (a stock is traded, air temperature is recorded, etc.) comes in discrete bits. It makes both mathematical and practical sense to similarly discretize the space in which the Ising model will live.

Mathematically, instead of having a random object defined for all points of the $d$-dimensional space $\mathbb{R}^d$, it will be defined only on the vertices of the $d$-dimensional cubic lattice $\Lambda$, like in Figure 3.

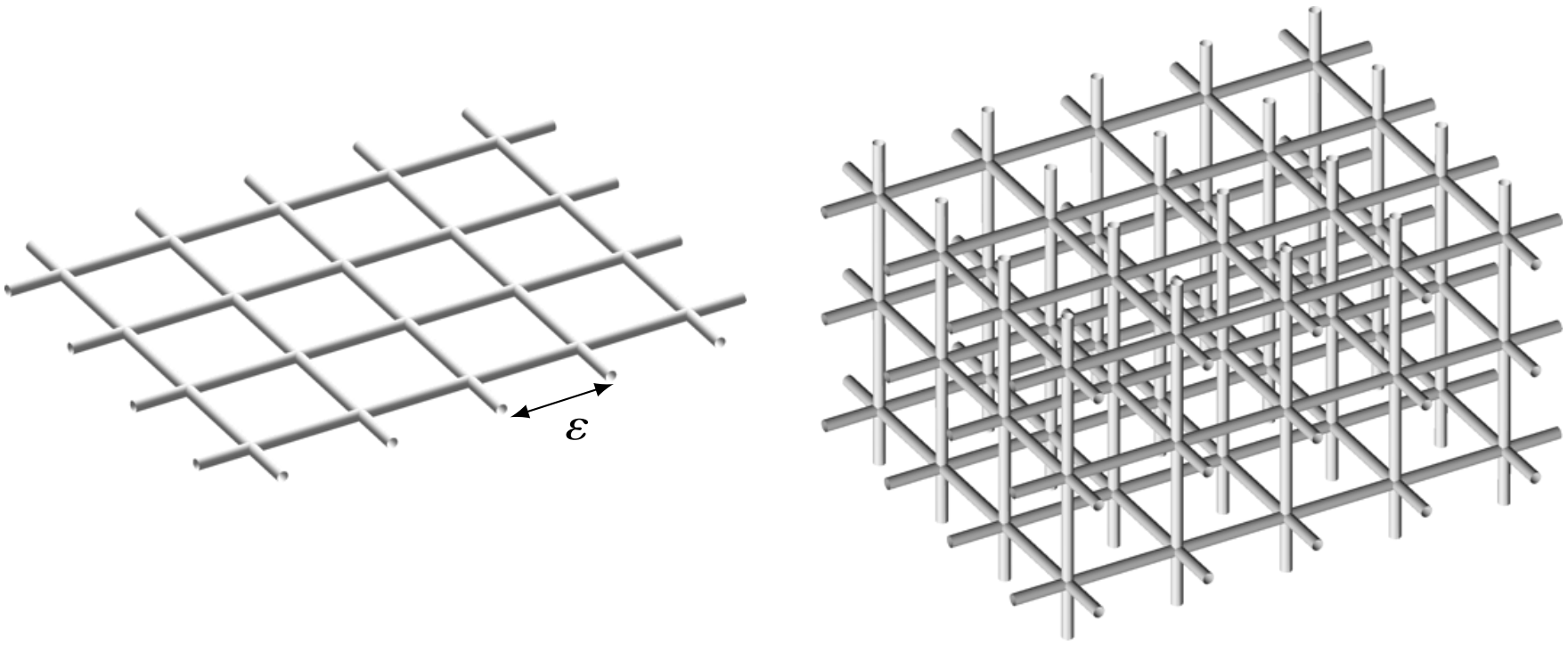

(3)

7 It still takes very exciting mathematics and lots of deep ideas to *prove* the behavior is Gaussian in dimensions $d \geq 4$.

Let $\varepsilon$ denote the mesh size of this lattice. Then the vertices of the lattice are the points whose coordinates are integer multiples of $\varepsilon$, that is,

$$\text{vertices}\,(\Lambda) = \varepsilon\mathbb{Z}^3 = \{(\varepsilon n_1, \varepsilon n_2, \varepsilon n_3)\} \subset \mathbb{R}^3\,, \tag{4}$$

where $n_1, n_2, n_3 \in \mathbb{Z}$ are integers. We will use the words "lattice vertex" and "lattice point" interchangeably.

From the human scale point of view, we should imagine $\varepsilon$ is vanishingly small, like the atomic scale. So, on the human scale, $\Lambda$ is very dense. But from the atomic scale point of view, we can take $\varepsilon = 1$. For an infinite lattice, both points of view are mathematically completely equivalent.

### 2.3. Signs on a lattice

Now it is time to assign some fluctuating degrees of freedom to the vertices of the lattice. In the Ising model, one makes the simplest possible *binary* choice. That is, at every vertex $\mathsf{v} \in \Lambda$, there is a random variable $\sigma(\mathsf{v})$ that can take two possible values. The reader may choose any name she or he likes for these values: black/white, blue/red, 0/1, $\pm 1$, et cetera. We will stick to the convention that

$$\sigma(\mathsf{v}) = \pm 1\,, \tag{5}$$

and we will call these variables *signs*. For historical reasons, they are normally called *spins*, which may be rather confusing for those familiar with spins and not familiar with the history of the Ising model. One advantage of (5) over 0/1 and other choices is that it stresses the *symmetry* between two possibilities. This symmetry is very important in the Ising model.

Minimalistically, a fragment of a configuration of the 2-dimensional Ising model may be represented like this:

$$\begin{matrix} + & - & + & - & + \\ + & + & + & + & - \\ - & - & + & + & - \\ - & - & - & - & + \end{matrix} \tag{6}$$

Ising, and his adviser Lenz, created the model as a model of ferromagnetism, and they imagined a miniature magnet at every site of a lattice pointing in one of two possible directions. For them, (6) represented something like this:

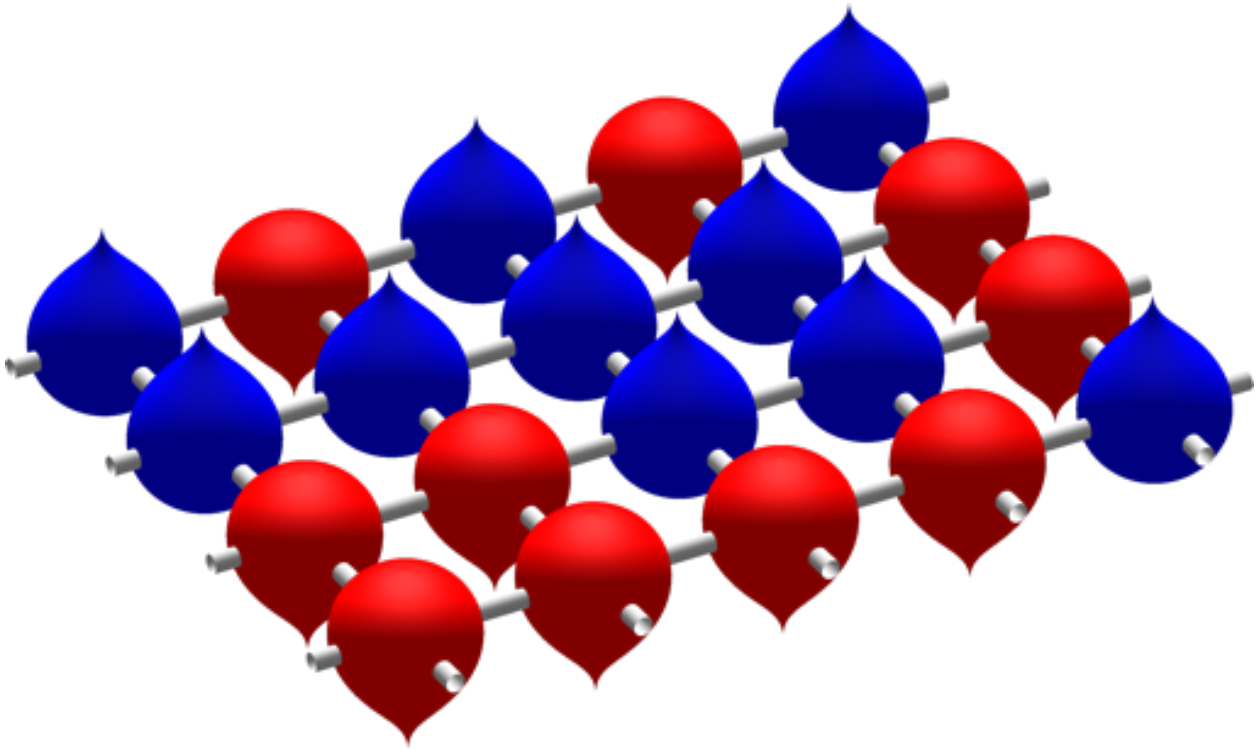

(7)

As it typically happens in mathematics, physics, and elsewhere, once introduced, mathematical models live their own life and follow their own logic. In particular, they may be used to understand phenomena that are very different from the originally envisioned applications.

History often preserves the context where people first stumble upon an important discovery. For instance, the mineral bauxite, the world's source of aluminium, is named after the village of Les Baux where it was first described by P. Berthier in 1821. In the Earth crust, it occurs mainly in places that are very far away from Provence and its current uses are probably very far from what Berthier could have envisioned.

### 2.4. Probabilities and Energy

A common misconception about probability theory is that all possible outcomes of a random event are equally probable. While this is a good approximation for coin tosses and dice rolls, this would not be a very interesting assignment of probabilities in the Ising model. Indeed, if all signs were equally likely, the picture in Figure (1) would be pure noise, indistinguishable from the noise introduced by the structure of the paper or the printing process. In particular, there would be no spatial structure to it, as it would be a collection of independent random bits, not affected by each other in any way. If they don't feel each other, they can be rearranged arbitrarily and hence there cannot be any significance to their particular spatial arrangment.

Ludwig Boltzmann and Willard Gibbs, the early architects of statistical physics, understood the connection between probabilities and *energy*. Energy is a central concept in physics which appears in the Newtonian, quantum, and statistical physics in slightly different, but compatible incarnations. Informally, it is supposed to be a universal equivalent that, just like ordinary human money, determines the intensity of any physical process. While people may have different attitudes towards money, energy is certainly making the physical world go round.

Without plunging into economic or metaphysical depths, mathematically energy is just a *function*

$$\{\text{configurations } \mathscr{C} \text{ of signs } \sigma\} \xrightarrow{\text{Energy}} \mathbb{R}\,, \tag{8}$$

that will be used to assign probabilities to configurations. Boltzmann and Gibbs understood that, in *equilibrium*, the probability of any configuration decays *exponentially* with its energy. To put these words into a formula, we have

$$\text{Prob}(\mathscr{C}) = \frac{1}{Z(T)} \exp\left(-\frac{\text{Energy}(\mathscr{C})}{T}\right). \tag{9}$$

In this formula, we have two proportionality coefficients $T$ and $Z(T)$ that both deserve a comment. Let's start with $T$.

Only dimensionless numbers make sense inside the exponential, but energy has physical dimension, namely

$$[\text{energy}] = [\text{mass}]\,[\text{length}]^2\,[\text{time}]^{-2}\,, \tag{10}$$

as exemplified by the familiar $\frac{1}{2}mv^2$ formula for the kinetic energy in the Newtonian mechanics. Therefore, we need a dimensional constant $T$ to convert energy into dimensionless numbers. The constant $T$ determines how fast probability decays with energy. Intuitively, it sets the scale of energy fluctuations. From

$$\frac{\mathrm{Prob}(\mathscr{C}_1)}{\mathrm{Prob}(\mathscr{C}_2)} = \exp\left(\frac{\mathrm{Energy}(\mathscr{C}_2) - \mathrm{Energy}(\mathscr{C}_1)}{T}\right) \tag{11}$$

we see that two configurations $\mathscr{C}_1$ and $\mathscr{C}_2$ make a comparable contribution only if the difference of their energies is not much larger than $T$.

From purely mathematical perspective, (9) may be taken as a *definition* of a statistical equilibrium, which depends on a constant $T \geq 0$ called the *temperature*[8].

The coefficient $Z(T)$ is defined so that the probabilities of all possible configurations $\mathscr{C}$ sum to 1. This is an interesting function of $T$ which, for historical reasons, is often called the *partition function*. Note that if we shift the energies of all configurations by the same constant $E_0$, then $Z(T)$ gets a factor of $e^{-E_0/T}$ and the probabilities do not change. In other words, only *energy differences* are important in (9). This is also clear from (11).

The case $T = 0$, interpreted using $e^{-\infty} = 0$, means the absolute zero temperature: only energy-minimizing configuration occur, and they are all equally likely.

### 2.5. Energy vs. Entropy

The dramatic plot of statistical mechanics is the competition between energy and *entropy*. Formula (9) gives preference to energy-saving configurations. They may each have a relatively large probability, but there are typically not so many of them. Having a close-to-minimal energy is a special property that most configurations will fail. But the competitive advantage of most configurations is that there are many of them.

To make a mathematical question out of this, we can ask how is the energy distributed in the system described by (9) ? The value of the energy in a random state of the system is a random variable, so it is a fair question. Anticipating the fact that in a system of large size the energy will scale linearly with a suitably defined volume $V$ of the system, it is better to look at energy $E$ per unit volume. One defines its entropy by

$$S(E) = \frac{1}{V} \ln \left(\text{number of states with } \frac{\mathrm{Energy}}{V} = E\right) . \tag{12}$$

---

**8** While (9) is a definition, it is still worth explaining why $T$ is called temperature. Imagine two systems in equilibrium at temperatures $T_1$ and $T_2$, respectively, which can exchange energy but otherwise don't interact. So, the configurations of the combined system are pairs $(\mathscr{C}_1, \mathscr{C}_2)$ and

$$\mathrm{Prob}((\mathscr{C}_1,\mathscr{C}_2)) = \mathrm{Prob}(\mathscr{C}_1)\ \mathrm{Prob}(\mathscr{C}_2)\,, \quad \mathrm{Energy}((\mathscr{C}_1,\mathscr{C}_2)) = \mathrm{Energy}(\mathscr{C}_1) + \mathrm{Energy}(\mathscr{C}_2)\,.$$

From (9) the combined system is in equilibrium if and only if $T_1 = T_2$. It thus suffices to check that (9) agrees with any other definition of a temperature for any one standard system, such as the ideal gas. Note that many thermometers work by putting some standard probe in contact and equilibrium with the system in question.

Here, again, we normalize the logarithm by the volume $V$ because we expect the counts of different possible states of the system to grow exponentially with the volume $V$. An equivalent of (12) is inscribed on Boltzmann's tombstone in Vienna's Zentralfriedhof.

Formula (12) is a definition, like formula (9). From these definitions, we conclude

$$\mathrm{Prob}\left(\frac{\mathrm{Energy}}{V} = E\right) = \frac{1}{Z(T)} \exp\left(\frac{V}{T} \underbrace{(-E + TS(E))}_{\text{maximize}}\right), \tag{13}$$

where we hid the inessential proportionality factor in gray.

When $V$ is very large, only those energies that minimize $E - TS(E)$, a quantity known as free energy, will be observed in the system. Not the ones that simply minimize $E$. The character of this minimum depends on the temperature. For $T = 0$, only the energy counts, and we get strict energy minima. For $T = \infty$, energy means nothing and entropy decides. For other values of $T$, both energy and entropy count, in different proportions. We will see this principle in action in the Ising model.

### 2.6. Interactions in the Ising model

Now it is time to specify the energy function in the Ising model. Let $\mathscr{C}$ be a configuration of signs. We can write it as a function

$$\sigma : \Lambda \longrightarrow \{\pm 1\}, \tag{14}$$

assigning each vertex $\mathsf{v} \in \Lambda$ a sign. When mathematicians talk about a function $\sigma$, they write $\sigma(\mathsf{v})$ to denote its value at the argument $\mathsf{v}$, and use the symbol $\sigma$ to denote the "whole" function. A configuration in the Ising model is a function (14) and we don't need another symbol $\mathscr{C}$ to denote it. What we need is to assign a number to it that will be called $\mathrm{Energy}(\sigma)$.

The spatial structure of the lattice will be taken into account by declaring that only neighboring signs interact. That is,

$$\mathrm{Energy}(\sigma) = \sum_{\text{edges } \mathsf{v}—\mathsf{v}'} E(\sigma(\mathsf{v}), \sigma(\mathsf{v}')) \tag{15}$$

where the edges are the edges in the lattice (3), the vertices $\mathsf{v}$ and $\mathsf{v}'$ are the two endpoints of a given edge, and $E(\pm 1, \pm 1)$ is some interaction energy of the neighboring spins to be specified momentarily.

Note that all edges contribute equally to (15), no matter where in the lattice they occur and in which of the coordinate directions they are pointing. In other words, the interactions in (15) are as homogeneous and as isotropic as the presence of a lattice in space allows.

It remains to specify 4 numbers $E(\pm 1, \pm 1)$. Since we want plus and minus to be symmetric, we need to have

$$E(1,1) = E(-1,-1), \quad E(1,-1) = E(-1,1),$$

where the latter equality also follows from the symmetry of the interaction between two neighbors. Recall that the overall shift of energy changes nothing and note that the overall scale of energy is equivalent to rescaling the temperature. In the end, there are no meaningful free parameters left, and we can set

$$E(\sigma(\mathsf{v}), \sigma(\mathsf{v}')) = -\sigma(\mathsf{v})\sigma(\mathsf{v}') \,. \tag{16}$$

The minus sign here means that signs lower the energy by being equal[9]. In other words, the function (14) likes to be a constant function, or the spins in (7) like to be pointing in the same direction. This desire, however, is not expressed globally, but only through local interaction of each sign with its immediate neighbors.

A careful reader may have been worried for a long time by the sum in (15) being an infinite sum of $\pm 1$'s for the infinite cubic lattice $\Lambda$. This worry is well justified and related to some core mathematical and physical issues. We will devote many pages below to dealing with it carefully. For now, let's replace the infinite lattice $\Lambda$ by any finite piece of it or any finite graph. Then (15) is a finite sum, the probabilities in (9) are well-defined, and we have defined the Ising model in finite volume.

### 2.7. Clusters and interfaces

Grouping together neighboring vertices of the same sign, we get *clusters* of pluses and minuses, as in the following figure:

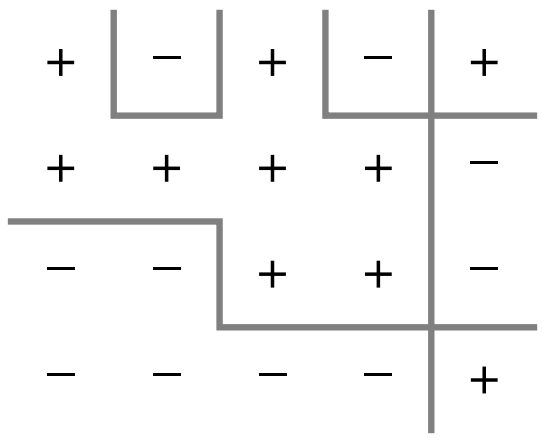

(17)

The boundary between the clusters is the *interface* between pluses and minuses. It is a $(d-1)$-dimensional object glued out of sides of a unit square/cube, so a path for $d = 2$, a surface for $d = 3$, et cetera. One component of the interface is highlighted in Figure (1). For $d = 3$, the interface may look something like Figure (48) in Section 3.3.3.

From (15) we have

$$\text{Energy}(\sigma) = \text{const} + 2\,\text{Area}(\text{interface}) \,, \tag{18}$$

where area (or length) is the $(d-1)$-dimensional lattice area, meaning that each side of the unit cube has area 1. This means that the Ising model can be interpreted as describing a fluctuating lattice interface, where the energy of the interface is its lattice area.

The lattice area has some peculiarities compared to usual area in $\mathbb{R}^d$. For instance, in $\mathbb{R}^2$, any two points are joined by a unique shortest path — a straight line segment, while the shortest path enclosing a given volume is a circle. For the 2-dimensional lattice distance,

---

9 The opposite choice of sign in (16) mathematically means negative temperature and, for other discretization of space, e.g. the triangular lattice, may correspond to a very different physics. For the cubic lattice, however, it can be reduced to the minus sign by flipping half of the signs in a checkerboard fashion.

there are many shortest paths connecting two points. Indeed, any path that goes up/right from the start to the finish in the following picture has a minimal length:

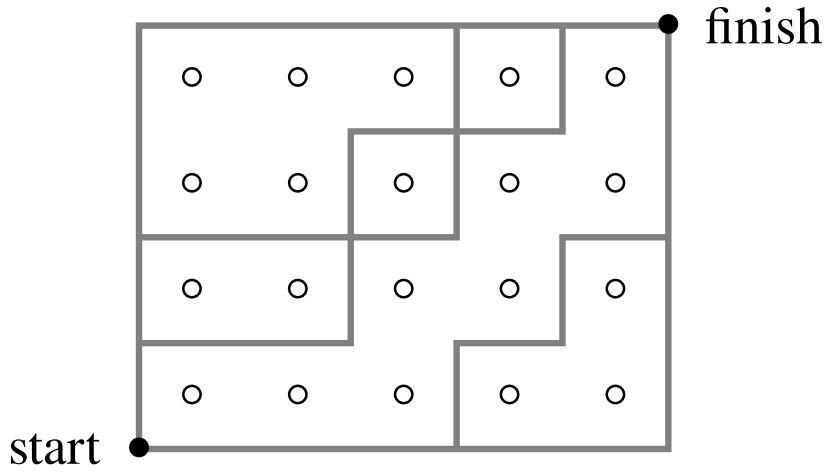

(19)

The shortest lattice paths enclosing a given volume are close to squares, not circles. Similar features persist in all dimensions.

In this narrative, we will be concerned with the *phase transition* in the Ising model that happens at a certain *critical temperature* $T_c$. These terms will be introduced properly below. For now we remark that for $T < T_c$, the Ising model reflects, due to the peculiarities of the lattice area, the behavior of materials that are similarly anisotropic, for instance crystals. Indeed, for a crystal, a square or a cube is the shape one should expect to see, not the sphere.

At $T = T_c$ however, this anisotropy disappears, a rather remarkable phenomenon. In fact, at the critical temperature, not only rotation invariance is restored, but some further symmetries appear. This is an incredibly interesting topic, but since it it would be a side trip for our story, we refer the reader to [18, 32] for details.

This concludes our brief discussion of the Ising model in finite volume. It is time to make sense of the energy and probabilities for the whole infinite lattice. This will be our task in Section 3.

In the process of defining the Ising model, there were choices, and we always made the simplest possible nontrivial symmetric choice. A reader may get the impression we defined a little mathematical toy, a basic wooden block set, which may be good for play but seriously oversimplifies the nature. What is the place of the Ising model in the broader landscape of statistical physics ? An interested reader will find an introductory discussion of this question in Appendix A. In short, mathematical physicists believe the Ising model provides a *universal* description of a very large class of phenomena in which a $\pm 1$ symmetry becomes *broken* below a certain temperature.

## 3. Gibbs measures

### 3.1. Definition

Our goal now is to define probabilities in the Ising model on the *infinite* cubic lattice. More precisely, we want to know what is the probability to see any particular pattern $\pi$ of signs in any given finite subset $\Omega$ of $\Lambda$. For instance, for $d = 2$, we want to know

$$\text{Prob}\left(\sigma\Big|_{\substack{\text{a fixed } 3\times 3\\ \text{square } \Omega}} = \boxed{\begin{matrix} + & - & + \\ + & + & - \\ - & - & + \end{matrix}}\right) = \,?\in [0,1]\,. \tag{20}$$

Mathematicians denote by $\sigma\big|_{\Omega}$ the restriction of a function $\sigma$ to a subset $\Omega$ of arguments. We will sometimes call the subset $\Omega$ a *window*. If we have a finite window into an infinite system, it is reasonable to ask what is the probability to see some pattern $\pi$ in it.

If $\Omega \subset \Omega'$, then the probabilities for the smaller window $\Omega$ are determined from the probabilities for the larger window $\Omega'$. Therefore, it is enough to define the probabilities for larger and larger cubes

$$\Omega_L = [-L, \dots, L]^d \subset \mathbb{Z}^d = \Lambda\,, \quad L = 1, 2, 3, \dots\,, \tag{21}$$

because any finite subset of $\Lambda$ is contained in some $\Omega_L$. For $d = 2$, the square $\Omega_1$ looks like the square in Figure (20).

The main issue with formula (9) for the infinite lattice was that the energy (15) is infinite. Recall, however, that the important thing in physics is not the energy itself, but rather the *difference* in energies and note that energy difference

$$\Delta\,\text{Energy} = \text{Energy}(\sigma) - \text{Energy}(\sigma') \tag{22}$$

is well defined if $\sigma$ and $\sigma'$ differ only at finitely many vertices.

Let's look at the example in (23), where the difference in signs is circled:

$$\sigma_+ = \begin{matrix} & \vdots & \vdots & \vdots & \\ \cdots & + & - & + & \cdots \\ \cdots & + & \oplus & - & \cdots \\ \cdots & - & - & + & \cdots \\ & \vdots & \vdots & \vdots & \end{matrix}\,, \qquad \sigma_- = \begin{matrix} & \vdots & \vdots & \vdots & \\ \cdots & + & - & + & \cdots \\ \cdots & + & \ominus & - & \cdots \\ \cdots & - & - & + & \cdots \\ & \vdots & \vdots & \vdots & \end{matrix}\,, \tag{23}$$

Assuming this is the only difference, that is, assuming that the dots in (23) represent some choice of signs for $\sigma_+$ and an *identical* choice for $\sigma_-$, we can compute the energy difference as follows:

$$\text{Energy}(\sigma_+) - \text{Energy}(\sigma_-) = 4\,. \tag{24}$$

Indeed, only the edges incident to the circled vertices change their energy, and their energy is $3 - 1 = 2$ for $\sigma_+$ and $1 - 3 = -2$ for $\sigma_-$. In exactly the same way, we can determine the change of energy if we flip some signs in the *interior* of the window $\Omega_L$ for any $L$.

Denote

$$\pi_+ = \boxed{\begin{matrix} + & - & + \\ + & + & - \\ - & - & + \end{matrix}}\,, \qquad \pi_- = \boxed{\begin{matrix} + & - & + \\ + & - & - \\ - & - & + \end{matrix}}\,. \tag{25}$$

We may interpret formula (9) as saying that

$$\frac{\text{Prob}(\sigma\big|_{\Omega_1} = \pi_-)}{\text{Prob}(\sigma\big|_{\Omega_1} = \pi_+)} = \exp\left(4T^{-1}\right) . \tag{26}$$

More generally, if $\pi$ and $\pi'$ are signs pattern in $\Omega_L$ that differ only in the interior, we may interpret formula (9) as saying that

$$e^{\text{Energy}(\pi)/T}\,\text{Prob}(\sigma\big|_{\Omega_L} = \pi) = e^{\text{Energy}(\pi')/T}\,\text{Prob}(\sigma\big|_{\Omega_L} = \pi') . \tag{27}$$

Note this formula operates with finite quantities. Moreover, (27) are *linear* equations on the probabilities.

Here comes the important moment. We say that an assignment of probabilities to the events that $\sigma\big|_{\Omega_L} = \pi$ is a *Gibbs measure*[10] if it satisfies (27) for all $L$ and all $\pi$ and $\pi'$ that differ only in the interior of $\Omega_L$. This key mathematical definition goes back to 1960s and the work of Roland Dobrushin, Oscar Lanford, and David Ruelle. And, no, Gibbs measures were not studied by Gibbs.

Note the change in perspective. Instead of saying that the formula (9) gives the probabilities, we have rewritten (9) as a system of equations that the probabilities have to satisfy. As with any equations, one naturally wonders: do they have a solution ? If they do, how many solutions are there ?

The existence of Gibbs measures for any temperature $T$ is a very general and soft mathematical fact, see Section 4.4.4. The question of how many Gibbs measures are there for a given value of $T$ is really the central question for us in these notes.

Formula (9) was meant to describe a statistical system in equilibrium at temperature $T$. But an infinite system can be in *many different* equilibria at given $T$, unable to fluctuate from one to another due to an infinite energetic cost. Concretely in the Ising model, each sign $\sigma(\mathsf{v})$ likes to be the same as its neighbors. Hence, a strong local preference for $+1$ or $-1$ may be self-reproducing in fluctuations. It could be a preference for either $+1$ or $-1$, and if there such a preference, then the system is stuck with it. The reader will probably have no difficulty thinking of real-life examples of this phenomenon.

Anticipating the fact that there may be many Gibbs measures at a given temperature, we will denote by $\boldsymbol{\mu}$ a Gibbs measure and write $\boldsymbol{\mu}(A)$ for the probability that $\boldsymbol{\mu}$ assigns to some event denoted by $A$. For example, $A$ can say that $\sigma\big|_{\Omega_L} = \pi$.

In practice, it is convenient to use the averages

$$\langle\sigma(\mathsf{v}_1)\sigma(\mathsf{v}_2)\ldots\sigma(\mathsf{v}_n)\rangle_{\mu} = \boldsymbol{\mu}(\text{this product equals } 1) - \boldsymbol{\mu}(\text{it equals } -1) \tag{28}$$

with respect to $\boldsymbol{\mu}$, where $\mathsf{v}_1, \ldots, \mathsf{v}_n$ are some vertices of $\Lambda$. The averages (28) are called *correlation functions*, and when one wants to stress the number of different lattice points involved, one talks about $n$-point correlation functions.

---

10 The word measure denotes a very important concept in mathematics, which we will leave without a proper discussion. The power of measure theory lies in being able to measure (meaning, assign some version of length, volume, probability, etc.) rather general sets. In our case, the probability is assigned to simple events of the form $\sigma\big|_{\Omega_L} = \pi$ and we hope the reader will have no difficulty thinking about this.

In general, the averages (also known as expectations, or integrals) with respect to a Gibbs measure $\boldsymbol{\mu}$ are defined as follows. Let $f(\sigma)$ be a function that depends on finitely many signs $\sigma(\mathsf{v}_i)$, $\mathsf{v}_i \in \Lambda$. Then $f$ takes finitely many values $f_j$, and we can define

$$\langle f \rangle = \sum_j f_j \, \boldsymbol{\mu}(f = f_j) \,. \tag{29}$$

In measure theory, general integrals with respect to a measure $\boldsymbol{\mu}$ are defined by approximating the integrand $f$ by functions taking finitely many values.

### 3.2. High temperature

#### 3.2.1.

It is easiest to start the discussion of Gibbs measures at the infinite temperature $T = \infty$. Since $1/T = 0$, energy disappears from (27) and we conclude that all sign patterns are equally likely. In other words, each sign is an independent symmetric coin toss. This is the complete description of the unique Gibbs measure for $T = \infty$.

#### 3.2.2.

For high enough temperatures, the unique Gibbs measure can be written as a series in the inverse temperature

$$\beta = \frac{1}{T} \,, \tag{30}$$

following a very general *perturbation theory* ideas, used everywhere in mathematical physics.

#### 3.2.3.

Let $\boldsymbol{\mu}_0$ be a Gibbs measure at inverse temperature $\beta_0$, from which we want to construct a Gibbs measure $\boldsymbol{\mu}$ at inverse temperature $\beta \approx \beta_0$. Let us first consider a finite piece $\Omega \subset \Lambda$ of the infinite lattice. For a finite graph $\Omega$, the unique Gibbs measure at inverse temperature $\beta$ is defined by (9). We can transform this definition as follows

$$\langle f(\sigma) \rangle_{\Omega,\beta} = \frac{\sum_\sigma e^{-\beta \,\mathrm{Energy}(\sigma)} f(\sigma)}{\sum_\sigma e^{-\beta \,\mathrm{Energy}(\sigma)}} \tag{31}$$

$$= \frac{\sum_\sigma e^{-\beta_0 \,\mathrm{Energy}(\sigma)} e^{(\beta_0-\beta)\,\mathrm{Energy}(\sigma)} f(\sigma)}{\sum_\sigma e^{-\beta_0 \,\mathrm{Energy}(\sigma)} e^{(\beta_0-\beta)\,\mathrm{Energy}(\sigma)}} \tag{32}$$

$$= \frac{\left\langle e^{(\beta_0-\beta)\,\mathrm{Energy}} f(\sigma) \right\rangle_{\Omega,\beta_0}}{\left\langle e^{(\beta_0-\beta)\,\mathrm{Energy}} \right\rangle_{\Omega,\beta_0}} \,, \tag{33}$$

where the summation in (31) and (32) ranges over all possible values of signs $\sigma(\mathbf{v})$ for $\mathbf{v} \in \Omega$.

Since the number $\Delta\beta = \beta - \beta_0$ is small, it may be useful to expand the exponentials in (33) in a series, using

$$e^x = 1 + x + \frac{x^2}{2} + \cdots + \frac{x^n}{n!} + \dots \,. \tag{34}$$

**3.2.4.**

For the infinite lattice $\Lambda$, formula (33) will seemingly run into the old problem of energy being infinite for an infinite lattice. However, it may happen that the infinities cancel between the numerator and denominator in (33) in each term of the series expansion in powers of $\Delta\beta$.

For concreteness, let's examine the first order of the expansion of a 1-point correlation function $\langle\sigma(\mathsf{v}_1)\rangle_\mu$. We have

$$e^{(\beta_0-\beta)\,\text{Energy}} = 1 + \Delta\beta \sum_{\text{edges } \mathsf{v}_2 — \mathsf{v}_3} \sigma(\mathsf{v}_2)\sigma(\mathsf{v}_3) + \dots . \tag{35}$$

where dots stand for terms of degree 2 or more in $\Delta\beta$. Therefore,

$$\left\langle e^{(\beta_0-\beta)\,\text{Energy}}\sigma(\mathsf{v}_1)\right\rangle_{\mu_0} = \langle\sigma(\mathsf{v}_1)\rangle_{\mu_0} + \Delta\beta \sum_{\text{edges } \mathsf{v}_2 — \mathsf{v}_3} \langle\sigma(\mathsf{v}_1)\sigma(\mathsf{v}_2)\sigma(\mathsf{v}_3)\rangle + \dots . \tag{36}$$

Dividing (36) by the average of (35), we obtain

$$\begin{aligned}\langle\sigma(\mathsf{v}_1)\rangle_\mu &= \langle\sigma(\mathsf{v}_1)\rangle_{\mu_0} + \\ &+ \Delta\beta \sum_{\text{edges } \mathsf{v}_2 — \mathsf{v}_3} \Big(\langle\sigma(\mathsf{v}_1)\sigma(\mathsf{v}_2)\sigma(\mathsf{v}_3)\rangle_{\mu_0} - \langle\sigma(\mathsf{v}_1)\rangle_{\mu_0}\langle\sigma(\mathsf{v}_2)\sigma(\mathsf{v}_3)\rangle_{\mu_0}\Big) + \dots\end{aligned} \tag{37}$$

The sum over the edges in (37) is infinite. However, if the edge $\mathsf{v}_2 — \mathsf{v}_3$ is far away from the vertex $\mathsf{v}_1$, we expect the corresponding signs to be approximately independent random variables. Approximate independence means that

$$\langle\sigma(\mathsf{v}_1)\sigma(\mathsf{v}_2)\sigma(\mathsf{v}_3)\rangle_{\mu_0} \approx \langle\sigma(\mathsf{v}_1)\rangle_{\mu_0}\langle\sigma(\mathsf{v}_2)\sigma(\mathsf{v}_3)\rangle_{\mu_0} . \tag{38}$$

So, the difference in (37) measures how close $\sigma(\mathsf{v}_1)$ and $\sigma(\mathsf{v}_2)\sigma(\mathsf{v}_3)$ are to being independent or, equivalently, how much they are *correlated*. If they decorrelate sufficiently fast with the distance between $\mathsf{v}_1$ and $\mathsf{v}_2$ then the sum in (37) will be convergent. Similar considerations apply to all other higher terms in the expansion (37).

**3.2.5.**

Given two random variables $f_1$ and $f_2$, the difference

$$\langle f_1 | f_2\rangle = \langle f_1 f_2\rangle - \langle f_1\rangle\ \langle f_2\rangle \tag{39}$$

is called their covariance. For example, the covariance of $f_1 = \sigma(\mathsf{v}_1)$ and $f_2 = \sigma(\mathsf{v}_2)\sigma(\mathsf{v}_3)$ with respect to $\mu_0$ appears in (37).

In statistical physics, it is typical for covariance (39) to decay if $f_1$ and $f_2$ depend on spatially separated arguments, like in (37). If this decay is exponential in the spatial separation then we will say that the model has an exponential decay of correlations or is exponentially decorrelated.

### 3.2.6.

There are higher analogs of the covariance, involving three of more arguments. For instance, one defines

$$\langle f_1|f_2|f_3\rangle = \langle f_1 f_2 f_3\rangle - \langle f_1 f_2\rangle\langle f_3\rangle - \langle f_1 f_3\rangle\langle f_2\rangle - \langle f_2 f_3\rangle\langle f_1\rangle + 2\langle f_1\rangle\langle f_2\rangle\langle f_3\rangle\,. \quad (40)$$

These are called *cumulants* and are related to the combinatorial principle of inclusion–exclusion. They measure finer mutual dependences between 3 or more random variables and appear naturally in perturbation series for the following reason.

The general formula for cumulants may be obtained from the identity

$$\ln\left\langle e^{f_1+f_2+\dots}\right\rangle = \sum_n \frac{1}{n!}\sum_{i_1,\dots,i_n}\left\langle f_{i_1}|f_{i_2}|\dots|f_{i_n}\right\rangle\,, \quad (41)$$

in which one expands the exponential as in (34) and equates terms containing the same functions $f_i$. In particular, let us replace $f_1$ in (41) by $tf_1$ and compute the value of $\frac{\partial}{\partial t}$ at $t=0$. We get

$$\frac{\left\langle f_1\, e^{f_2+f_3+\dots}\right\rangle}{\left\langle e^{f_2+f_3+\dots}\right\rangle} = \sum_n \frac{1}{n!}\sum_{i_1,\dots,i_n\geq 2}\left\langle f_1|f_{i_1}|f_{i_2}|\dots|f_{i_n}\right\rangle\,. \quad (42)$$

Voila, this is just what we need in (33), with $f_1 = \prod \sigma(\mathsf{v}_i)$ and

$$f_2 + f_3 + \dots = \sum_{\text{edges } \mathsf{v}—\mathsf{v}'} \sigma(\mathsf{v})\sigma(\mathsf{v}')\,.$$

Later in Section 4.3 we will meet random variables for which all cumulants with $n \geq 3$ vanish, meaning that that any $\langle f_1 \dots f_n\rangle$ may be written entirely in terms of the expectations $\langle f_i\rangle$ and the covariances $\left\langle f_i|f_j\right\rangle$. Such random variables are called *Gaussian*. See Section A.3.4 for more on this. In a certain precise technical sense, nonzero cumulants with $n \geq 3$ measure the nonlinearity of the model.

### 3.2.7.

Going back to the special case $\beta_0 = 0$ and the unique Gibbs measure $\mu_0$ at $T = \infty$, we observe that signs at different lattice sites are totally independent for $\mu_0$. Thus the $\Delta\beta$ term in (38) is simply zero. In fact, great simplifications happen for $\mu_0$ and a nice convergent combinatorial series can be written down for $\mu$ provided the inverse temperature $\beta$ is sufficiently small[11]. In this high-temperature range, the Gibbs measure remains unique.

Since $\mu_0$ and the energy are invariant under flipping all signs, this property is inherited by the perturbation series. The invariance of $\mu$ under flipping all signs also follows from

---

**11** Since the energy (35) is a product of terms like

$$e^{\beta\sigma(\mathsf{v}_2)\sigma(\mathsf{v}_3)} = \cosh(\beta)(1 + \tanh(\beta)\sigma(\mathsf{v}_2)\sigma(\mathsf{v}_3))\,,$$

it is more convenient to write this series in powers of the hyperbolic tangent of $\beta$:

$$\tanh(\beta) = \beta - \frac{1}{3}\beta^3 + \frac{2}{15}\beta^5 - \dots\,.$$

its uniqueness. It follows that

$$\big\langle \text{product of odd number of } \sigma(\mathsf{v}_i) \big\rangle_{\text{high T}} = 0\,, \tag{43}$$

and in particular that

$$\big\langle \sigma(\mathsf{v}) \big\rangle_{\text{high T}} = 0\,, \tag{44}$$

for any $\mathsf{v}$. The expected value of a single sign in (44) is the simplest measures of a possible $\pm 1$ asymmetry of a Gibbs measure. It is a very important parameter of the Gibbs measure called *magnetization*.

Also note the uniqueness of the high-temperature Gibbs measure implies it is invariant under *shifts* of the lattice $\Lambda$. This *translational invariance* is an important property for a Gibbs measure to have or not to have. In a translation-invariant Gibbs measure, the magnetizations at all vertices of the lattice are equal.

#### 3.2.8.

One should stress that the perturbative series expansion for (33) has no guarantee of success in general. In particular, it will fail when either $\boldsymbol{\mu}_0$ or $\boldsymbol{\mu}$ have long-range correlations, meaning that the signs at distant vertices do not become decorrelated sufficiently fast. Needless to say, these are precisely the situations of maximal interest and significance !

### 3.3. Low temperature

#### 3.3.1.

What about the opposite case $T = 0$ ? Equation (27) has the following meaning at $T = 0$

$$\text{Energy}(\pi) > \text{Energy}(\pi') \quad \Rightarrow \quad \text{Prob}(\sigma\big|_{\Omega_L} = \pi) = 0\,. \tag{45}$$

In other words, if the energy of a configuration can be lowered by flipping finitely many signs, then its probability vanishes.

In terms of the interface between the pluses and minuses, formula (18) says that it should be *minimal*, meaning that its length/area cannot be made any smaller by finite modifications[12].

What minimal interfaces can we think of ? First, there is the empty interface. If the interface is empty then all signs are equal. A measure $\boldsymbol{\mu}_+$ that assigns probability 1 to the configuration in which all signs are $+$, and zero probability to all other configurations, is a Gibbs measure at $T = 0$. No randomness is a special case of randomness and it may happen that the probability of just one particular configuration equals 1.

Since with an empty interface all signs can be $+$ or all signs can be $-$, we have *two* Gibbs measures $\boldsymbol{\mu}_\pm$ already. But there is more. For instance, the plane $x_1 = \frac{1}{2}$, or any plane

**12** For $d = 1$, there is a difference between finite modifications of signs and of the interface. We will consider finite modifications of the interface.

of the form $x_i = \frac{1}{2}$ + integer, $i = 1, 2, 3$, defines a minimal interface, see Figure (46)

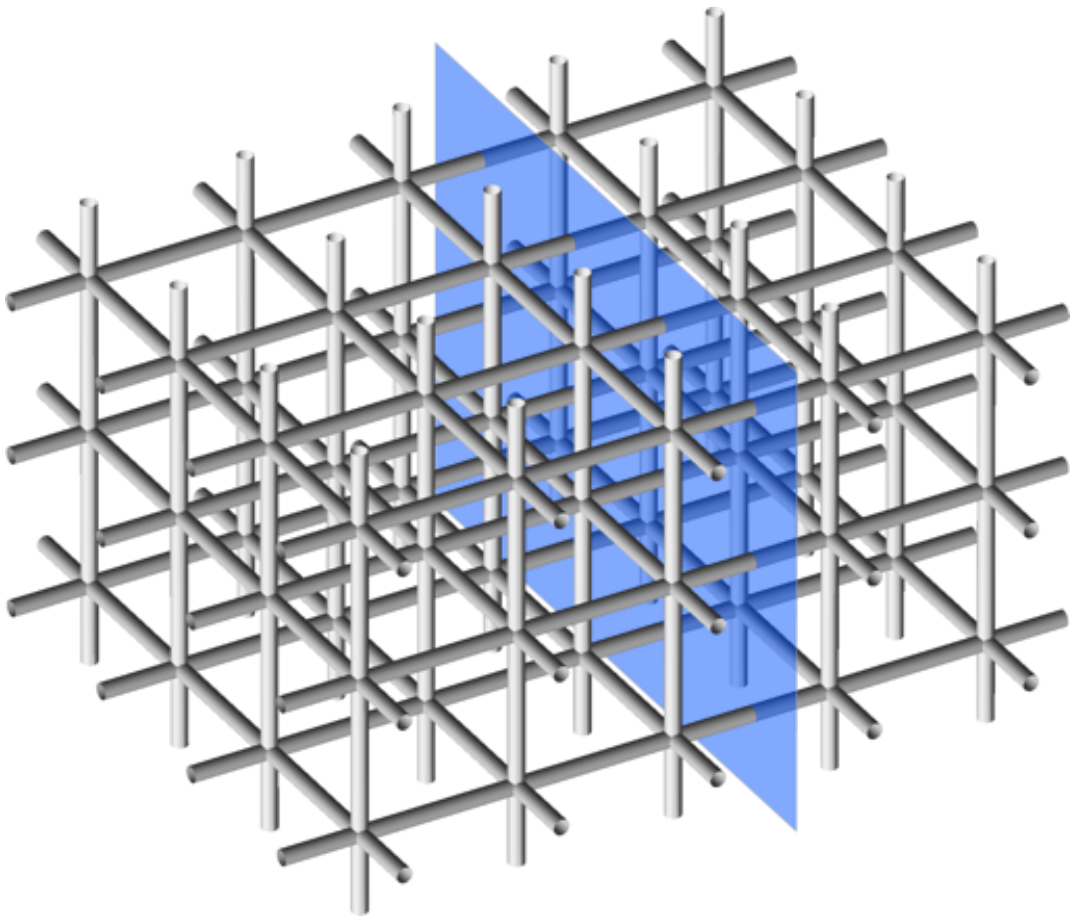

(46)

We can put pluses on either side of such wall and this gives many more zero-temperature Gibbs measures, all of which we will denote by $\boldsymbol{\mu}_{\text{wall}}$. They are not translation-invariant and, in fact, they can be all taken to one another by a symmetry of the lattice $\Lambda$.

A curious reader may think about more Gibbs measures at $T = 0$, but it is already clear that there are plenty of them. They very visibly break the symmetries between $\pm 1$ and also between different lattice points.

Recall that at $T = \infty$ we have a total *disorder*, which persisted to all high temperatures and manifested itself, in particular, by the vanishing magnetization. By contrast, the $T = 0$ measures exhibit a very strong spatial *order*.

**3.3.2.**

When there is more than one Gibbs measure, the following point should be kept in mind. Let $\boldsymbol{\mu}_1$ and $\boldsymbol{\mu}_2$ be two Gibbs measures. Then their mixture of the schematic form

$$\boldsymbol{\mu}_{\text{mix}} = 0.71\boldsymbol{\mu}_1 + 0.29\boldsymbol{\mu}_2 \,, \tag{47}$$

where 0.71 can be replaced by any number between 0 and 1, is also a Gibbs measure. Indeed, it assigns probabilities in $[0, 1]$ to all events and satisfies the linear equations (27) from the definition of a Gibbs measure.

What does the equation (47) mean ? Imagine there are two different labs in a physics department, labeled by $i = 1, 2$. In the lab number $i$, our physical system is kept in a state described by the measure $\mu_i$. We go to a *random* lab and do the measurement. If our chance to go to the first lab is 0.71 then the outcome of our measurement will be described by (47).

Clearly, this is rather silly. It is quite unnatural and adds nothing to our understanding of the system. Therefore, when there is more than one Gibbs measure, people usually restrict their attention to those Gibbs measures that cannot be nontrivially written in the form (47). They are called *extremal* or *pure*.

**3.3.3.**

How will the $T = 0$ Gibbs measure perturb for small positive $T$ ? The interface between plus and minus no longer has to be minimal, but every time its area increases by 1 the probability decreases by $e^{-2/T}$. It is therefore natural to organize the expansion in powers of $e^{-2/T}$ which is a small parameter for $T$ positive and small.

For example, we may perturb the wall in Figure (46) as as follows

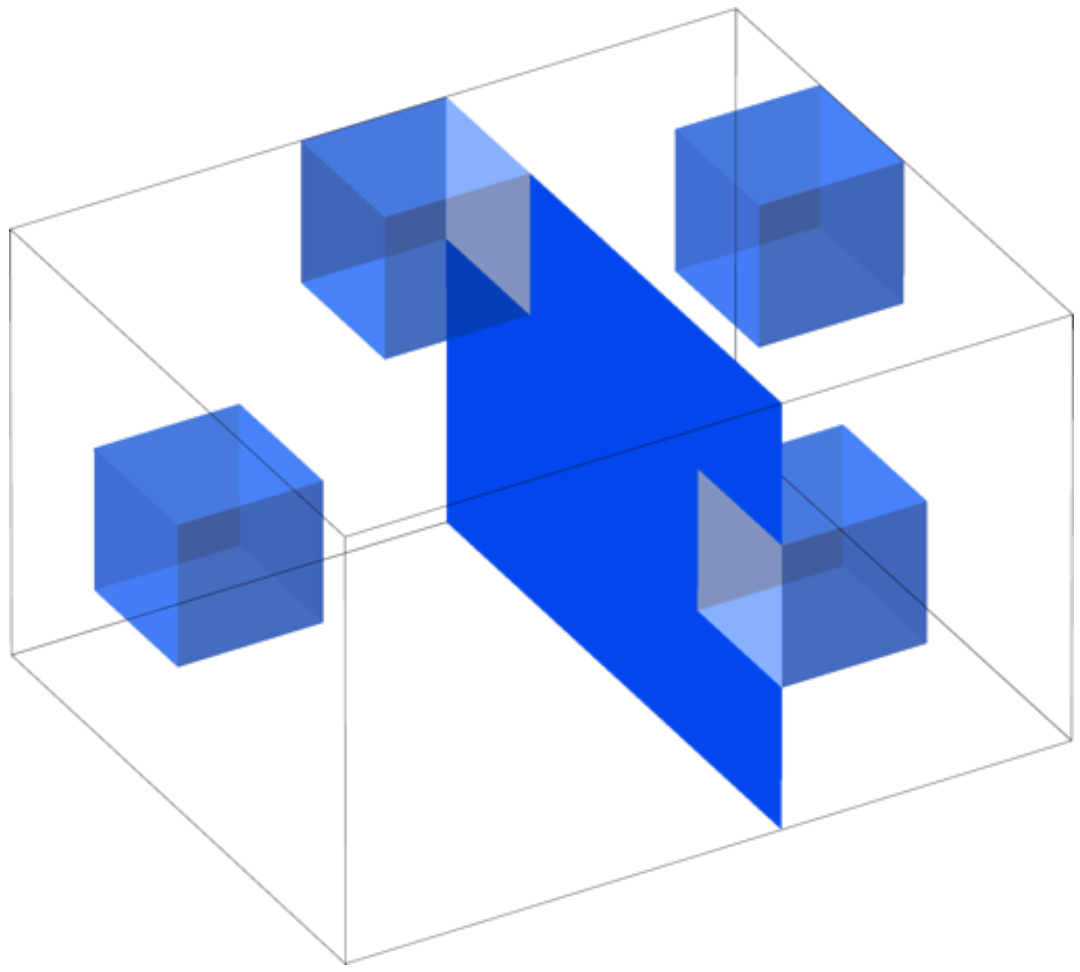

(48)

This perturbation increases the area by 20, so it counts with the weight $\left(e^{-2/T}\right)^{20} = e^{-40/T}$.

The fate of the perturbation series for zero temperature Gibbs measures is *different* in different dimensions.

**3.3.4.**

For $d = 1$, there are only the measures $\mu_\pm$ at $T = 0$. Their perturbation series breaks down at the very first term. Indeed, let $\mu_+$ denote the hypothetical Gibbs measure for which the signs are positive at the positive infinity of the lattice $\mathbb{Z}$.

The $e^{-2/T}$ term in the perturbation series for $\mu_+$ is then a sum over all configurations like this

. . . 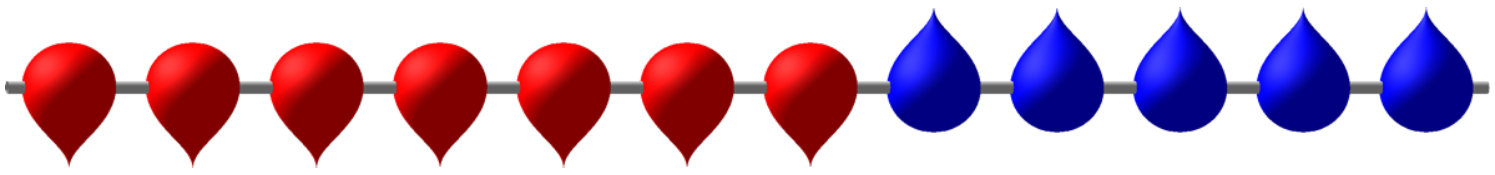 . . . , (49)

where the switch of signs can occur at any place. This gives infinitely many equal terms that affect the sign at any lattice point.

In fact, for $d = 1$, the high temperature disordered behaviour happens for all $T > 0$. No wonder we cannot access any positive temperature by a perturbation of the $T = 0$ description. Historically, Ernst Ising studied precisely the $d = 1$ case and reached this conclusion in his 1924 dissertation.

The absence of order for any $T > 0$ in $d = 1$ lead to a certain temporary dip of interest in the Ising model. See [17] for a much more informative account of the many chapters of the Ising model's history.

### 3.3.5.

For $d \geq 2$, the perturbation series for $\boldsymbol{\mu}_{\pm}$ converges ! This was first noted, in essense, in 1936 by Rudolf Peierls. Peierls observed that the number of relevant interfaces of given area is bounded by $C^{\text{Area}}$ for some constant $C$. This makes the series converge as long as $e^{2/T} > C$ and proves that, for $d \geq 2$, the Ising model can exhibit both order and disorder, depending on the temperature.

### 3.3.6.

For $d = 2$, the perturbation series for $\boldsymbol{\mu}_{\text{wall}}$ breaks down at the first possible term when the length is allowed to increase by 2. We already discussed in Section 2.7 that the lattice length has just too many minimizers, and this is another consequence of this fact. In fact, in $d = 2$, the measures $\boldsymbol{\mu}_{\pm}$ can be shown to be the *only* pure Gibbs measures for $T > 0$.

By constrast, for $d \geq 3$, the series for $\boldsymbol{\mu}_{\text{wall}}$ converges. The corresponding measures were first studied by Dobrushin and bear his name [10, 11]. This means that at low temperature and in dimensions $d \geq 3$, the Ising model can *break* both the $\pm$ symmetry *and* the symmetries of the lattice. One says that some symmetry $\mathbf{g}$ of the system is broken by a Gibbs measure $\boldsymbol{\mu}$ if $\mathbf{g}$ takes $\boldsymbol{\mu}$ to another Gibbs measure, different from $\boldsymbol{\mu}$.

## 3.4. Critical temperature

We have talked about the behaviour of the Ising model at high and low temperatures, respectively. This behaviour differs strikingly. At high temperatures, we have a homogenous disorder. The system expresses no $\pm 1$ preference and looks the same everywhere. At low temperatures, vertices prefer one sign over the other and this preference may change from vertex to vertex.

What happens for temperatures in the middle ? This question must be on everybody's mind by now. Is there some intermediate range of temperatures for which yet another qualitatively different behavior is observed ?

For the Ising model, and related models of statistical physics, there is exactly one *critical* value $T_c$ of the temperature at which the balance between order and disorder, energy and entropy tips. A numerical simulation, done by Stanislav Smirnov, may help visualize this transition.

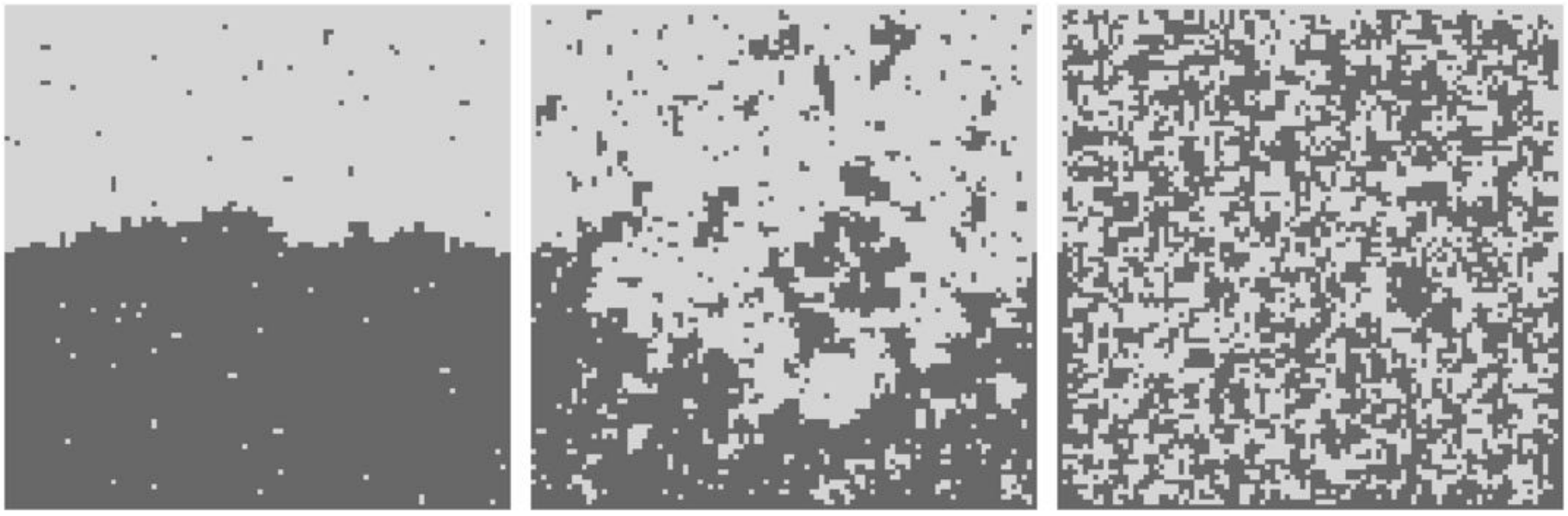

(50)

In (50) we see the $d = 2$ Ising model simulated on the $100 \times 100$ grid for $T < T_c$, $T = T_c$, and $T > T_c$ respectively. For a finite piece of the lattice, the notion of *boundary conditions* is important. Formula (27) tells us about the probabilities of different signs patterns inside the square, while signs along the boundary may be in principle assigned arbitrarily. In (50) and also in (1), we have +1 along one the lower half of the boundary and −1 along the upper half, see Figure (51), as if we are trying to simulate the nonexistent Gibbs measure $\mu_{\text{wall}}$ in $d = 2$.

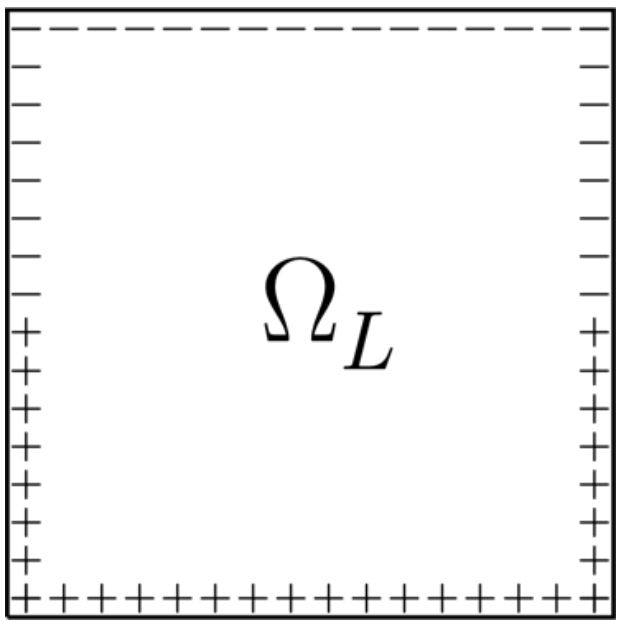

(51)

Such boundary condition allows one to focus on the properties of the interface by singling out one special component of the interface that goes from one vertical side of the boundary to another. This interface is clearly visible for $T < T_c$, has essentially evaporated for $T > T_c$, and it is both visible and very strongly fluctuating at the critical point $T = T_c$.

All features of the $T > T_c$ are microscopic, everything happens on the lattice scale. The $T < T_c$ has one macroscopic feature — the interface, but its fluctuations are again small in size and would not be visible from far away[13]. What is common between the $T \neq T_c$ pictures is that the signs at different lattice points become independent *exponentially* fast with the lattice distance. The rate of this exponential decay sets the typical scale of observed features in both pictures.

By contrast, the $T = T_c$ picture has some features on *all scales*, enabled by the slow *polynomial* decorrelations. In fact, the lattice mesh $\varepsilon \to 0$ limit of the critical $d = 2$ Ising model is invariant not just under scaling, but under all *conformal* transformations. These are transformations that look like scaling and rotation in a very small neighborhood of any point[14]. For the proof of conformal invariance of the $d = 2$ critical Ising model, Stanislav Smirnov was awarded the 2010 Fields Medal, the only Fields Medal previously awarded for the study of the Ising model. While our focus in these notes is on $d = 3$, we hope a curious reader will open [18, 32] for more on conformal invariance. After the groundbreaking 1984 work of three Alexanders: Belavin, Polyakov, and Zamolodchikov, conformal invariance and the language of Conformal Field Theory (CFT) grew to be the most powerful tool for understanding 2-dimensional critical phenomena.

---

**13** It may be useful to explain, in terms of Figure (50), why there is no Gibbs measure $\mu_{\text{wall}}$ in $d = 2$ for $0 < T < T_c$. If we fix any finite window at exactly the middle height, the interface will pass over it or under it with probability almost $\frac{1}{2}$. As a result, we will observe $\mu_\pm$ in our window with probabilities $\frac{1}{2}$ as the square in (51) grows to infinity. For $d = 3$ and $T < T_c$, the interfaces fluctuates less and we will see it in the window as the measure $\mu_{\text{wall}}$..

**14** It is easier to define conformal transformation as the transformations that preserve all angles between curves. They can scale and rotate by different amounts in the vicinity of different points. The case $d = 2$ is special in that there is an abundance of such transformations.

It is the job of a statistical physicist to predict macroscopic properties of materials from the macroscopically invisible fluctuations that take place on the atomic scale. It is a very, very interesting job, with its challenges and rewards, and the study of the Ising model at $T \neq T_c$ is no exception. But the statistical physicist's finest hour is when she or he gets to describe a system that does exhibit macroscopic fluctuations. Phase transitions are these kind of phenomena. In particular, the order/disorder phase transition in the Ising model certainly packs more excitement, and is much more widely applicable, than what happens at $T \neq T_c$. So, what happens at $T = T_c$ ? This question calls for the start of a new section.

## 4. What happens at $T = T_c$ ?

### 4.1. Critical Gibbs measures

With our focus on Gibbs measures in this narrative, it is clear what our next question is going to be. Is there one or are there many Gibbs measures at $T = T_c$ ? This question may be phrased as *continuity* of the phase transition. Indeed, if there are many Gibbs measures at $T = T_c$, there will be one of them, say, $\mu_c$, which is *not* the $T \downarrow T_c$ limit of the unique high-temperature measure $\mu_{\text{high T}}$. Thus, for a system in state $\mu_c$, a tiny increase in temperature will lead to a jump to $\mu_{\text{high T}}$, meaning a jump in physical properties.

In a live or online class, it may be a good idea to take a poll on this question. Do you think it is going to be continuous ? Or not ? Phase transitions come in both flavors in nature. When the water melts or boils, its properties change discontinuously. At the pressure of 1 atmosphere, water boils at 100°C. Increasing the pressure increases the boiling point monotonically until, at the pressure of 217.7 atmospheres, we reach a very special point called the water critical point. After it, the difference between liquid and vapor disappears. When going through this point, the properties of the system remain continous. Admittedly, this is a much more delicate example than simply boiling the water.

Another example of a continuous phase transition is the Curie critical point, the original motivation for the introduction of the Ising model. Magnets loose their magnetic properties when heated. For an iron magnet, this happens at 770°C, and the loss of magnetic properties is continuous. The Ising model is not a particularly convincing model of magnetism for several reasons, so we should be careful with drawing conclusions from this example.

### 4.2. The Potts model

A very important difference between a magnet and Ising model signs is that magnetization is a vector that can be rotated in all possible ways. These rotational symmetries are very different from the simple $\pm 1$ symmetry of the binary degrees of freedom in the Ising model. Rotations form a continuous Lie group. Importantly, rotations can be arbitrarily small.

Closer to the Ising model are the models with a larger, but still finite symmetry group. The most important example is the $Q$-state Potts model, the Ising model being the

$Q = 2$ case of the Potts model. In the Potts model, the function

$$\sigma : \Lambda \longrightarrow \{1, 2, \ldots, Q\}, \tag{52}$$

can take $Q$ possible values, and the energy has the same form (15) with

$$E(a, b) = \begin{cases} E_=, & a = b, \\ E_{\neq}, & a \neq b. \end{cases} \tag{53}$$

This is invariant under all $Q!$ possible permutations of the values in (52). As with the Ising model, the exact values of constants in (53) are not important as long as $E_{\neq} > E_=$.

Simulations of the critical $Q$-state Potts model by V. Beffara for $d = 2$ and $Q = 2, 3, 4, 5, 6, 9$ may be seen in Figure (54).

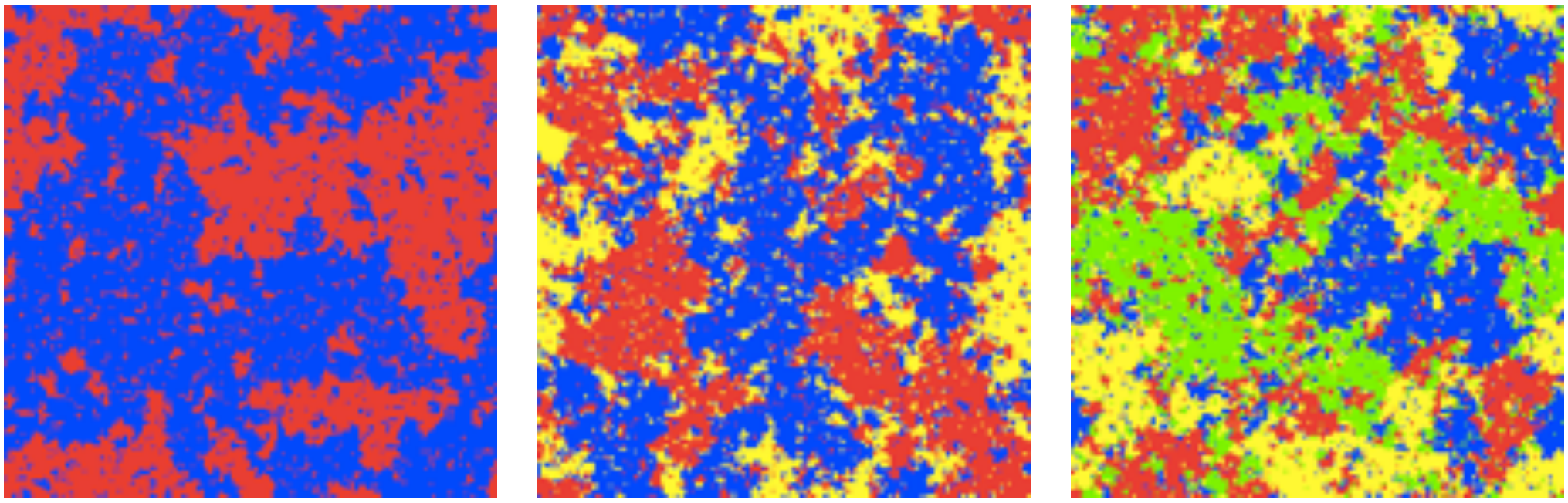

(54)

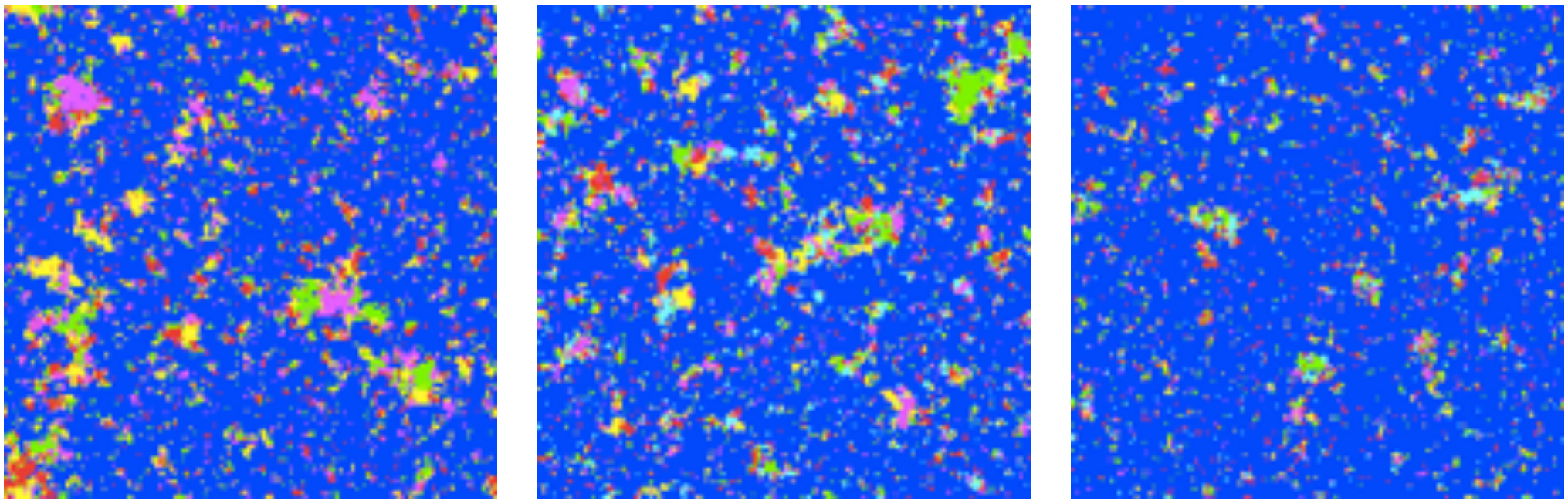

Different colors represent clusters of different values of $Q$. The reader will notice a clear difference in behaviour between the top $Q \leq 4$ row and the bottom $Q \geq 5$ row.

### 4.3. Theorems

Our goal on the preceding pages was to ignite the reader's interest in what happens in the 3-dimensional Ising model at the critical temperature. To add to the suspense, we start with the following fundamental result of Hugo Duminil-Copin and his collaborators in the $d = 2$ case.

**Theorem 1** ([19,20]). *The phase transition in the $Q$-state Potts model for $d = 2$ is continuous if and only if $Q \leq 4$.*

Theorem 1 is a logical conjunction of two different results proven using two different sets of tools. The continuity for $Q \leq 4$ is proven in the paper [19] by Hugo Duminil-Copin, Vladas Sidoravicius and Vincent Tassion, building, in particular, on the earlier solo work [12] of Hugo Duminil-Copin. The discontinuity for $Q > 4$ is proven in the paper [20] by Hugo Duminil-Copin, Maxime Gagnebin, Matan Harel, Ioan Manolescu, and Vincent Tassion.

After the depth of the continuity question has been underscored yet another time by Theorem 1, we can finally state the following result of Michael Aizenman, Hugo Duminil-Copin, and Vladas Sidoravicius.

**Theorem 2** ([1]). *The phase transition in the $d = 3$ Ising model is continuous.*

It is expected that the phase transition for the $d = 3$ Potts model is discontinuous for $Q \geq 3$, see [16]. As to the higher dimension we have the following result. Recall we have met the Gaussian random fields in Section 3.2.5, see also Section A.3.4 in the Appendix. The scaling limit refers to the limit when we take correlations, suitably scaled, for larger and larger spatial separations. It describes what we would actually observe on our human scale[15].

**Theorem 3** ([2]). *The scaling limit of the critical Ising model in $d = 4$ is Gaussian.*

I hope the readers share the narrator's sense of awe at this absolutely amazing mathematics and join me in warmest congratulations on it being recognized by the Fields Medal. I also hope the readers got the sense that today's mathematics is not just extraordinarily powerful, but also concrete, understandable, and fun, once one finds the right idea and the right point of view. While finding that right point of view is not at all easy, my biggest hope is to have inspired my youngest readers to believe that mathematics can be beautiful and rewarding, both as a subject and as a profession. Maybe this is also a good place for me to thank Hugo Duminil-Copin, Stanislav Smirnov, and Martin Hairer for this special opportunity to be introduced to their wonderful subject.

## 4.4. Contours of proofs, seen in the distance

### 4.4.1.

We hope the reader agrees that the majestic view of Theorem 1, 2, and 3 was worth the uphill hike through the foothills of the Ising range. We also hope the reader will not be discouraged to learn that a much longer and steeper climb is needed to get a good view of the actual mathematics that goes into the proof of these theorems. As we stressed at the beginning, having a mathematical proof is a measure of our understanding of the model and, certainly, understanding is a great reward for any effort.

To help the reader master the subject, there are brilliant expositions available, in particular by Hugo Duminil-Copin himself. In [16], the reader will find a very fun, colorful, and engaging explanation of Theorem 1, 2 and many other results.

15 if we lived in corresponding number of spatial dimensions

### 4.4.2.

Let us start with Theorem 2 and a discussion of the basic logic of how something like this could be proven. One logical point we should make from the very beginning is that we *do not* know the value of $T_c$ for $d = 3$.

It is a gift of nature to mathematical physicists that many fascinating and highly nontrivial exact results can be obtained in the Ising, Potts, and related models when $d = 2$. In particular, it is known that for the square lattice Potts model we have

$$T_c = \frac{2}{\ln(1+\sqrt{Q})}\,, \quad d = 2\,. \tag{55}$$

This goes back to the 1941 work of Kramers and Wannier in the $Q = 2$ Ising case and is proven by Vincent Beffara and Hugo Duminil-Copin in [5] in general.

Formulas like this are extremely sensitive to the exact lattice formulation of the model, and other $d = 2$ models presumably converging to the same critical CFT at their (unknown!) critical point loose the magic. In addition to being a huge help in the study of the Ising and the Potts models proper, exact results very much contributed to how mathematical physicists think about their subject in general. We will say a few words about them below.

Nothing of the kind was ever discovered for $d = 3$, and there are many different strong hints that the physics, and the mathematics, in the plane and in the space are just different.

### 4.4.3.

Recall our discussion of the $T = 0$ Gibbs measures and note that, of all possible Gibbs measures, the measure $\boldsymbol{\mu}_+$ clearly has the most pluses, while the measure $\boldsymbol{\mu}_-$ has the least possible number of them. This basic comparison persists to all temperatures. All possible Gibbs measures are, in a certain precise mathematical sense, sandwiched between $\boldsymbol{\mu}_-$ and $\boldsymbol{\mu}_+$. Hence, the continuity question may be phrased as

$$\boldsymbol{\mu}_+ \stackrel{?}{=} \boldsymbol{\mu}_-\,, \quad T = T_c\,. \tag{56}$$

To see whether $\boldsymbol{\mu}_+ \stackrel{?}{=} \boldsymbol{\mu}_-$, one doesn't need to compute all correlation functions. Well-developed techniques in the subject reduce the question to the comparison of 1-point correlation functions, that is, magnetizations, at all vertices. Since $\boldsymbol{\mu}_\pm$ are both translation-invariant and differ by exactly the flip of all signs, the continuity question is equivalent to

$$\big\langle \sigma(\text{any one point}) \big\rangle_{\boldsymbol{\mu}_+} \stackrel{?}{=} 0\,, \quad \text{at } T = T_c\,. \tag{57}$$

This may sounds like we made good progress until we remind ourselves that we don't know the value of $T_c$, or any equation that determines this number. An approximate value of $T_c$ is known from numerical experiments, but it is not useful for us now. The only thing we know about $T_c$ is that

$$T_c = \inf \Big\{ T,\ \text{such that } \big\langle \sigma(\mathsf{v}) \big\rangle_{\boldsymbol{\mu}_+} = 0 \Big\}\,. \tag{58}$$

But since the $T \downarrow T_c$ continuity is precisely the crux of the matter, we didn't progress much. We will be just going in circles until we can relate the question (57) to something which is

either:

(A) true for ALL temperatures, or

(C) is manifestly CONTINOUS as $T \downarrow T_c$ . (59)

#### 4.4.4.

In the (C) category in (59), one can actually describe the limit of the unique high-temperature Gibbs measure as $T \downarrow T_c$. One has

$$\lim_{T \downarrow T_c} \mu_{\text{high T}} = \mu_{\text{free}} \Big|_{T=T_c} , \tag{60}$$

where $\mu_{\text{free}}$ is the *free boundary* Gibbs measure that can be constructed as follows.

There is a universal way to produce Gibbs measures for any temperature. Recall the discussion of the boundary condition from Section 3.4. For all sufficiently large $L$, say $L \geq 100$, fix some configurations of signs

$$\sigma\Big|_{\partial\Omega_L} = \tau_L , \quad L = 100, 101, 102, \ldots , \tag{61}$$

along the boundary $\partial\Omega_L$ of the cube $\Omega_L$ from (21) as in Figure (62). The values of $\tau_L$ can be all $+1$, can be all $-1$, can be like the Dobrushin's boundary conditions from (51), can be anything at all. In particular, they don't have to be related to each other for different values of $L$.

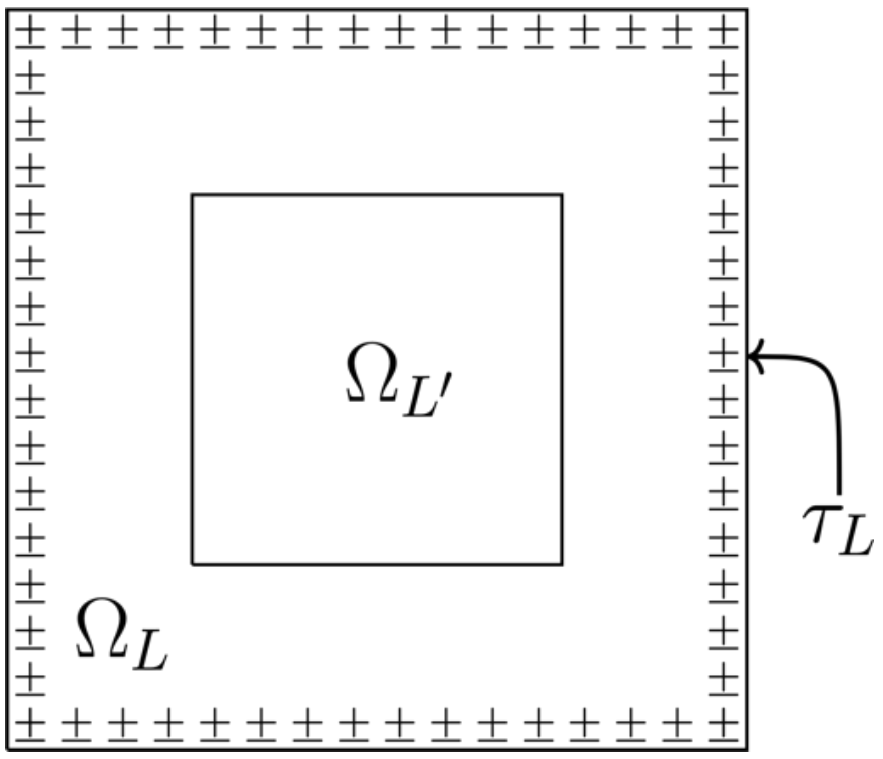

(62)

For all $L$, equation (27) defines a unique probability distribution $\mu_L$ for signs in $\Omega_L$ and thus probability distributions for signs in any smaller cube $\Omega_{L'} \subset \Omega_L$. Possible values of $\sigma\big|_{\Omega_{L'}}$ form a finite set and probability distributions on a finite (or compact) set are compact[16]. Hence there is a subsequence of $L$ for which the limit

$$\mu(\sigma|_{\Omega_{L'}} = \pi) = \lim_{L\to\infty} \mu_L(\sigma|_{\Omega_{L'}} = \pi) \tag{63}$$

exists for all $L'$ and $\pi$. Since it is a limit of solutions of (27), it is a Gibbs measure. Incidentally, this proves that the set of Gibbs measures is nonempty for any $T$.

---

16 The key property of a compact set that we need here is that any infinite sequence of elements of a compact set has a converging subsequence.

Any Gibbs measure may be obtained in this way. Indeed, if we make $\tau_L$ random and let it be distributed according to some Gibbs measure $\mu$ then $\mu_L = \mu$.

The measure $\mu_{\text{free}}$ is obtained when instead of fixing the boundary signs, we take all possible sign configuration in $\Omega_L$, weighted according to the temperature and their energy. This means we really sum over all possible configurations of $\pm$ along the boundary in (62) with the corresponding weights. The measure $\mu_{\text{free}}$ is not pure. In fact, it is known that [6]

$$\mu_{\text{free}} = \tfrac{1}{2}\mu_+ + \tfrac{1}{2}\mu_- \,, \qquad T < T_c \,. \tag{64}$$

While $\mu_\pm$ are at the two extremes of the set of the Gibbs measure, the measure $\mu_{\text{free}}$ is its very center. In particular, it is $\pm$-symmetric. It thus makes sense, and can be shown formally, that it is the $T \downarrow T_c$ limit of the unique high-temperature measure in (60).

### 4.4.5.

In the (A) category in (59), we would like a comparison between $\mu_+$ and $\mu_{\text{free}}$ which is valid for *all* temperatures.

From (64), we see that it is hopeless to compare magnetizations, as they will definitely differ below $T_c$. However, since $\mu_+$ and $\mu_-$ differ by a sign flip, their $n$-point correlation functions are equal for any *even* number $n$. This means that both above and below $T_c$ we have

$$\left\langle \prod_{i=1}^{n} \sigma(\mathsf{v}_i) \right\rangle_{\mu_+} = \left\langle \prod_{i=1}^{n} \sigma(\mathsf{v}_i) \right\rangle_{\mu_{\text{free}}} , \qquad T \neq T_c \,, \qquad n \text{ is } \underline{\text{even}} \,, \tag{65}$$

and thus it is a reasonable hope to extend this to $T = T_c$.

In fact, if (65) can be extended to $T = T_c$ for $n = 2$ that would be the end of the proof because of the following argument.

On the one hand, as a very special case of a general FKG inequality published in 1971 by Cees Fortuin, Pieter Kasteleyn, and Jean Ginibre [23], one has

$$\langle \sigma(\mathsf{v}) \rangle_\mu \langle \sigma(\mathsf{v}') \rangle_\mu \leq \langle \sigma(\mathsf{v}) \sigma(\mathsf{v}') \rangle_\mu \tag{66}$$

for any Gibbs measure $\mu$. For a translation-invariant measure $\mu$, it follows that[17]

$$\langle \sigma(\mathsf{v}) \rangle_\mu^2 \leq \lim_{\|\mathsf{v}-\mathsf{v}'\| \to \infty} \langle \sigma(\mathsf{v}) \sigma(\mathsf{v}') \rangle_\mu \,. \tag{67}$$

For $\mu_{\text{free}}$, the right-hand side of (67) can be seen to vanish at $T_c$ as a consequence of the $T \downarrow T_c$ continuity. If the 2-point functions for $\mu_{\text{free}}$ and $\mu_+$ are the same then (67) implies $\langle \sigma(\mathsf{v}) \rangle_{\mu_+} = 0$ at $T = T_c$ and we are done.

### 4.4.6.

What do we remember about the sign configurations if we forget all $n$-point correlations for $n$ odd ? It is easy to see that we remember precisely the clusters of equal signs

**17** As a side remark, the inequality in (67) is, in fact, an equality for $\mu_+$ and this is how Onsager's formula (71) for magnetization for $d = 2$ was originally derived.

which we talked about in Section 2.7, see Figure (68)

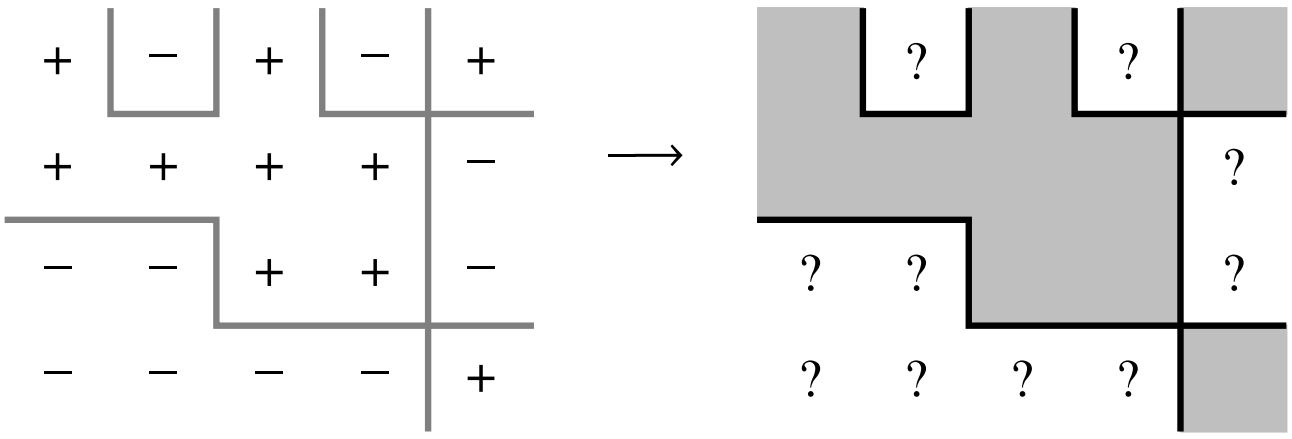

(68)

The Ising model thus becomes a *random cluster model*, the random object in which is a random partition of lattice vertices into clusters[18]. Such random cluster models play a very important role in mathematical physics and are closely related to various *percolation* models, see [16]. In a percolation model, the edges of a lattice, or of a more general graph, are kept or erased with some probabilities and the connected pieces of what remains are called the percolation clusters.

The analog of a magnetization for a random cluster model is the probability that two *neighboring* vertices $\mathsf{v}$—$\mathsf{v}'$ belong to the same cluster. A closely related quantity is $\langle \sigma(\mathsf{v})\sigma(\mathsf{v}')\rangle_{\mu}$, where $\mathsf{v}$—$\mathsf{v}'$ is an edge of the lattice. By an analysis reminiscent of how (56) is deduces from (57), the authors of Theorem 2 show:

$$\begin{array}{c} \langle \sigma(\mathsf{v})\sigma(\mathsf{v}')\rangle_{\mu_+} = \langle \sigma(\mathsf{v})\sigma(\mathsf{v}')\rangle_{\mu_{\text{free}}} \\ \text{for an edge } \mathsf{v}\text{—}\mathsf{v}' \end{array} \quad \Rightarrow \quad \begin{array}{l} \text{the random cluster models} \\ \text{for } \boldsymbol{\mu}_+ \text{ and } \boldsymbol{\mu}_{\text{free}} \text{ are equal}\,. \end{array} \tag{69}$$

Recall that the equality of the random cluster models implies the equality (65) for all $T$.

### 4.4.7.

It “only” remains to show that

$$\langle \sigma(\mathsf{v})\sigma(\mathsf{v}')\rangle_{\mu_+} - \langle \sigma(\mathsf{v})\sigma(\mathsf{v}')\rangle_{\mu_{\text{free}}} = 0\,. \tag{70}$$

for one edge $\mathsf{v}$—$\mathsf{v}'$. And this is where the *real* ascent or perhaps even flight begins and our excursion wraps up.

We will just say that the authors of [1] estimate the left-hand side in (70) using a certain auxiliary percolation model, the edges in which are kept or erased by a procedure that takes its input from the Ising model or, more precisely, from the *random current* representation of the Ising model. This random current representation may be compared and contrasted with the high-temperature expansion from Section 3.2.

Recall how we talked in the beginning about a mathematician’s freedom do introduce and use any auxiliary mathematical structure that may shed a new light on the question at hand. Just like a geometer is free to introduce any auxiliary construct, a mathematical physicist is free to introduce any auxiliary model, limited only by one’s own imagination. While there is no physical percolation happening in the Ising model, one can learn a great deal about the Ising model from the percolation model studied in [1].

---

**18** In specialized literature, the term *random cluster model* often refers to a particular class of models that are related to the Ising clusters by a further random refinement, see [16].

### 4.4.8.

The proof of Theorem 1 is very different and is based on certain highly nontrivial *exact* results for the square lattice $d = 2$ $Q$-state Potts model. From the early days of the Ising model to the modern heights of Theorem 1, exact results played a very important role in the development of statistical physics, and mathematical physics in general.

For instance, the continuity in the $d = 2$ Ising model case follows at once from the celebrated formula of Onsager (see e.g. [3, 4] for a historical account)

$$\langle \sigma(\mathrm{v}) \rangle_{\mu_+} = \left(1 - \frac{1}{\sinh^4(2/T)}\right)^{\frac{1}{8}}, \quad T \leq T_c\,, \tag{71}$$

where the formula (55) for the critical temperature is obtained solving the equation

$$\sinh(2/T_c) = 1\,.$$

The plot of this function can be seen in Figure (72)

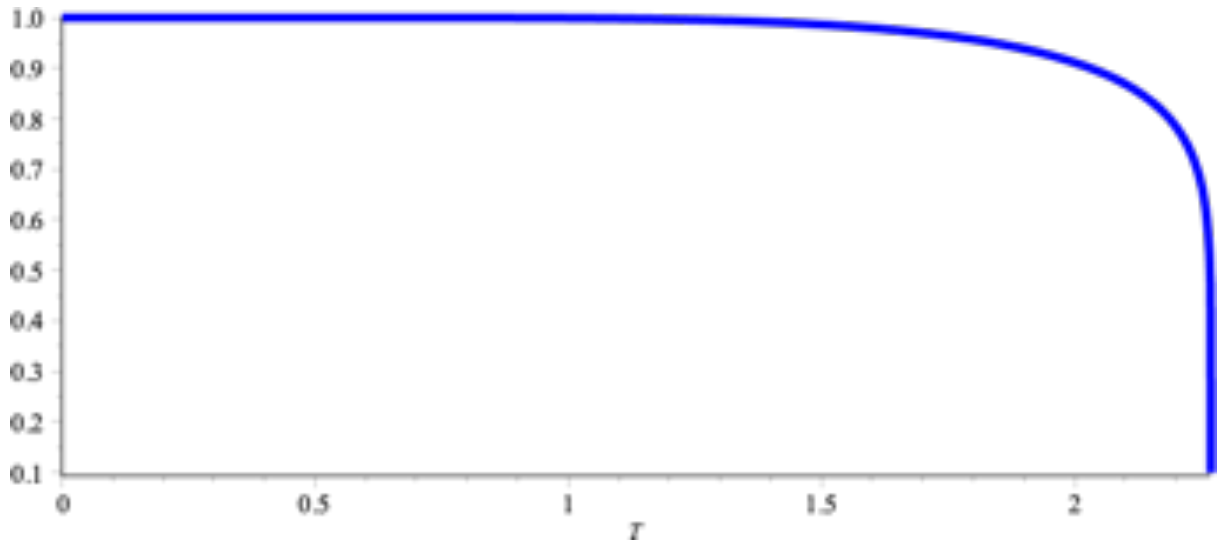

(72)

In addition to the continuity, we observe the remarkable fact that magnetization behaves like $(T_c - T)^{\frac{1}{8}}$, a result that found a deep explanation in conformal field theory [19].

### 4.4.9.

Ultimately, the algebraic structure responsible for exact computations is a certain infinite-dimensional symmetry algebra present in the $Q$-state Potts model on the square lattice. It extends to the discrete lattice level the infinite-dimensional symmetries of the CFT limit — a most remarkable phenomenon. The algebra in question is a $q$-deformation of the Lie algebra of $2 \times 2$ matrices with entries in Laurent polynomials in one variable. The parameter $q$ of this deformation is related to the parameter $Q$ by

$$q + q^{-1} = -\sqrt{Q}\,.$$

Hence the difference between $Q > 4$ and $Q \leq 4$ is the difference between $q$ being negative real and $q$ being a complex number on the unit circle $|q| = 1$.

While this is a strikingly beautiful story in mathematics, it doesn't really belong in our narrative, with our focus on the amazing $d = 3$ breakthrough achieved in a total absense

**19** Numerical bootstrap computations [21] predict that in $d = 3$ the magnetization behaves like $(T_c - T)^{\beta}$ where $\beta = 0.326419\ldots$. This is a much more robust and universal number than the critical temperature, but a good mathematical understanding of it awaits future generations of mathematical physicists.

of exact results. We do suggest, however, that the interested reader opens [8] for a general introduction to quantum group with a view towards their applications, [26, 29] for classic treatments by some of the key figures in the development of the subject, and maybe also [25] for a representation-theoretic take on the origin of the structures used in the $Q \leq 4$ part of the proof. All of this, of course, in addition to the brilliant exposition of the actual proofs in [16].

Asked whether he likes exact results or estimates better, Hugo Duminil-Copin says: "*I prefer estimates. They usually offer a more robust approach to critical phenomena, and I am as much as possible trying to obtain proofs that are not based on exact formulae.*"

"*I certainly agree with Hugo, inequalities are more versatile. But I am not sure we would have advanced so far without having some exact identities first. It is a miracle that there are any equalities concerning the Ising model. Onsager's calculation was shocking at the time, as it provided an exact formula for a function exhibiting a phase transition. This and later miraculous equalities, together with inequalities, forged our understanding of the Ising model. Hugo is a master of both inequalities and equalities, and has really moved the frontier of statistical physics with many beautiful theorems with equally beautiful proofs. I congratulate Hugo from all my heart, bravo!*", says Stanislav Smirnov.

"*Hugo Duminil-Copin's work has brought unprecedented clarity to our mathematical understanding of phase transitions in statistical mechanics. The elegance of his proofs truly makes them seem come straight out of The Book*", says Martin Hairer.

## 5. Further reading

Popular accounts of these and related developments include [9, 34]. See especially the popular piece [15] written by Hugo Duminil-Copin for the Oberwolfach's snapshots of modern mathematics.

We quoted many times from [16] and an interested reader is certainly advised to continue her or his exploration of the subject following these lectures. Among other survey articles written by Hugo Duminil-Copin one may list [13, 14, 17].

To anyone who can read French, Hugo Duminil-Copin wholeheartedly recommends the lectures [27] by Jean-François Le Gall and the book [33] by Wendelin Werner. Another very important book is the subject is [24] by Geoffrey Grimmett.

I hope the reader has a lot of fun studying these sources as well as the original articles including [1].

## A. The universal attraction of the Ising model

### A.1. Universality

It is hard to tell the atomic composition of a liquid by watching it evaporate or freeze. There is a good reason it took humans millenia to figure out the microscopic composition of macroscopic objects. Part of the reason is that a great many different microscopic systems have the same macroscopic behavior.

Molecules live on a nanometer (that is, $10^{-9}$m) scale and there is incredibly many of them in a macroscopic piece of any material ($18 = 2 + 16$ grams of $H_2O$ contain about $6 \cdot 10^{23}$, the Avogadro number, of molecules). It sounds completely impossible that their individual behaviour could be observed by us. Instead, we observe only the combined, or averaged, effect of myriads of molecules.

For instance, if we have a container of gas, we can measure the density, the pressure, the temperature, etc. These measure the average number of molecules[20] per unit of volume, the average force per unit area exerted by the gas on the wall of the container, and the average kinetic energy of molecules, respectively. In principle, we could measure more quantities, but the equation of state (an important concept in statistical physics) tells that temperature and pressure are enough. Add to this the vector of the wind, and there is no further information about the gas that a weather station can provide.

Mathematically, what does it mean that there is no further information ? Recall the concept of a Gibbs measure from Section 3. It assigns a probability to every event one can detect and hence an average, or expected value, to any observable quantity. If some finite number of these expectations already determine the whole Gibbs measure then they determine the outcome of *every* possible measurement in our system.

Going back to the gas, on a macroscopic scale, it is described by 2 scalars, temperature and pressure, and one vector, wind velocity. These vary in space and time if the gas is not in global equilibrium and are the ingredients in the mathematical models of motion of gases and fluids. In some sence, the job of a statistical physicist is to provide the arrow

$$\boxed{\text{description on the 1nm scale}} \xrightarrow{\substack{\text{statistical}\\ \text{physics}}} \boxed{\text{description on the 1m scale}}, \qquad (73)$$

connecting two very different kind of physics, and two communities of mathematical physicists studying the corresponding phenomena using very different models and mathematical tools. Since the source and the target in (73) are so very different, it is impossible for the target to be a faithful image of the source. To reiterate, we can't tell the atomic composition of air just by feeling a cool breeze. Put differently, an enormous amount of information is discarded by the arrow (73).

This loss of information is a win for a statistical physicist. It means there is no pressing need to study every possible scenario of microscopic interactions. People call it *universality*, meaning the macroscopic conclusions should hold universally and indepen-

20 If we have a mixture of several gases, there will be separate densities for each kind of molecules. These can be traded for the corresponding partial pressures.

dently of most microscopic details. Within each universality class, it is thus reasonable to restrict our attention to the simplest possible microscopic model.

Universality is a very important ingredient in how statistical physicists think about their subject. To be clear, it is always an enormous mathematical challenge to prove any universality statement rigorously. However, there is an appealing heuristic description of the universality classes based on the *renormalization group* idea of Kenneth Wilson. We will say a word about it below.

### A.2. Models like the Ising model

#### A.2.1.

How can stuff fluctuate in space ? We should have some fluctuating degrees of freedom, which we may describe by an $N$-tuple of numbers

$$\phi(x) = (\phi_1(x), \phi_2(x), \dots, \phi_N(x)) \in \mathbb{R}^N \,.$$

Here, the argument $x$ is a $d$-dimensional vector, which we will discretize to a lattice $\Lambda \subset \mathbb{R}^d$. It is good to visualize $\Lambda$ as a fine mesh $\varepsilon\mathbb{Z}^d \subset \mathbb{R}^d$ approximating the space $\mathbb{R}^d$ in the *continuous* limit $\varepsilon \to 0$.

In parallel with (14), we can write $\phi(x)$ as a random function

$$\Lambda \xrightarrow{\ \ \phi\ \ } \mathbf{\Phi} \subset \mathbb{R}^N \,, \tag{74}$$

where $\mathbf{\Phi}$ is the range of the possible values of $\phi$. In the Ising model, for instance, $N = 1$ and $\mathbf{\Phi} = \{\pm 1\}$.

#### A.2.2.

The interactions are described by an energy function, such as the energy (15) in the Ising model. In general, one imagines

$$\text{Energy} = \text{External potential} + \text{Pair potential} + \dots \,, \tag{75}$$

where dots stand for other possible interactions. We will assume (75) is translation-invariant. Then the first term has the form

$$\text{External potential} = \sum_{x\in\Lambda} U_1(\phi(x)) \,, \tag{76}$$

for some function $U_1$ on $\mathbf{\Phi}$ in (74). In the continuous $\varepsilon \to 0$ limit, the sum in (76) becomes the integral of $U_1(\phi)$.

In the Ising model, one can add such term. This is called Ising model in an *external field*. It breaks the $\pm 1$ symmetry and destroys the critical point. It is very interesting, however, to study the response of the critical Ising model to a small external field.

If (76) is the only nonzero term in (75), then from (9) we conclude that the values of $\phi(x)$ are independent identically distributed $N$-dimensional random variables with probability density function proportional to $e^{-U_1(\phi)/T}$. In space, this is a complete noise, with some nontrivial distributions of values, hence not something of great interest to us now.

### A.2.3.

A translation-invariant pair potential has the form

$$\text{Pair potential} = \sum_{x,y\in\Lambda} U_2(\phi(x), \phi(y), x-y)\,. \tag{77}$$

This term puts spatial interactions in (75). We may assume each term in (77) is $x \leftrightarrow y$ symmetric.

For the cubic lattice Ising model, $U_2$ vanishes unless $x - y = \pm e_i$, where $e_i$ are the coordinate vectors. In principle, one can allow next-nearest neighbors to interact, as well as lattice sites further away[21]. It is important, however, for the interaction to decay rapidly with the distance between $x$ and $y$. Models in which everything interacts with everything behave like a crowd and are usually well-described by the crowd average $\overline{\phi}$ fluctuating in some potential $U_{\text{effective}}(\overline{\phi})$ derived from $U_1$ and $U_2$.

For a pair of neighbors $\mathsf{v}$—$\mathsf{v}'$ in the Ising model, we can write

$$-\sigma(\mathsf{v})\sigma(\mathsf{v}') = -1 + \frac{1}{2}\frac{|\sigma(\mathsf{v}) - \sigma(\mathsf{v}')|^2}{\|\mathsf{v} - \mathsf{v}'\|^2} \tag{78}$$

because $\mathsf{v} - \mathsf{v}'$ is a unit vector and $\sigma$ takes values $\pm 1$. The fraction on the right in (78) is a lattice version of the square of the derivative of $\sigma$ in the direction of $\mathsf{v} - \mathsf{v}'$. Since an overall shift of energy does nothing, we see that the pair energy in the Ising model can be written as a discretization of $\frac{1}{2}\|\nabla\sigma\|^2$, where $\nabla$ denotes the gradient of the function.

In general, we may think of (74) as of discretization of a map $\mathbb{R}^d \longrightarrow \mathbf{\Phi}$, like in Figure (79).

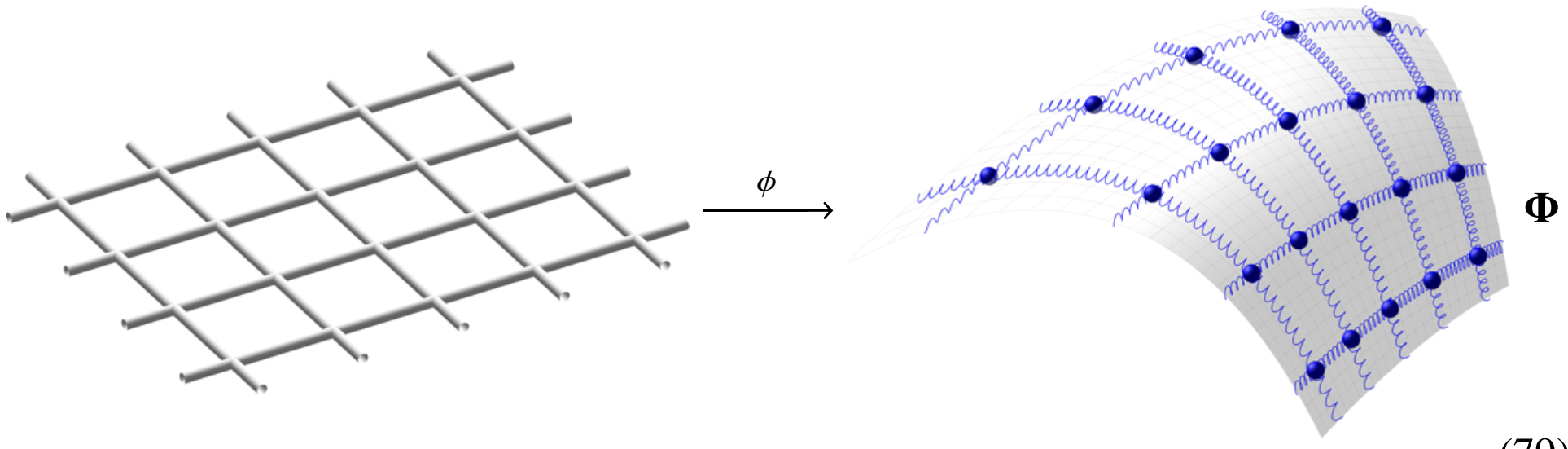

(79)

The role of the pair potential is to hold the values of this map together by putting an energy price on wild oscillations. One natural notion of energy for a continuous map is the Dirichlet energy $\frac{1}{2}\int \|\frac{\partial\phi}{\partial x}\|^2$, the construction of which, in general requires a metric in the domain and target of $\phi$. The anisotropic and anharmonic relatives of the Dirichlet energy are certainly possible and important in the description of materials with the corresponding properties.

**21** In fact, the authors of Theorem 2 prove it in much wider generality than described in these notes.

## A.3. Critical points

### A.3.1.

Wilson's idea[22] was that the arrow in (73) can be presented as a composition of many similar arrows that each change the scale by modest factor, such as 2

$$\begin{array}{ccccccccccc} 2^{-30}\, m & \longrightarrow & & & & & & & & & 1\, m \\ \| & & & & & & & & & & \| \\ 2^{-30}\, m & \longrightarrow & 2^{-29}\, m & \longrightarrow & 2^{-28}\, m & \longrightarrow & \dots & \longrightarrow & \dots & \longrightarrow & 1\, m\,. \end{array} \tag{80}$$

Stepping off the firm mathematical grounds for the rest of this section, we may hypothesize that the change of scale by 2 corresponds to some *renormalization* transformation

$$(U_1, U_2, U_3, \dots) \xrightarrow{\ R\ } (U_1', U_2', U_3', \dots)\,, \tag{81}$$

where $U_3$ corresponds to possible triple interactions in (75) etc. Since the arrow in (73) is the transformation $R$ raised to some very large power, we should put two theories in the same universality basket if they become identical aften many iterations of (81).

### A.3.2.

While heuristic, this argument underscores the importance of *scale-invariant* models. If there really was a well-defined transformation (81), such theories will be its fixed points. It makes sense that the result in (73) is scale invariant, since we certainly expect the same macroscopic description to be valid at both the $1m$ and $2m$ scales.

The invariance here should be understood up to redefinition of the fields. Indeed, if $\phi_1(x)$ is measured in meters then $R$ should act on it by $\phi_1(x) \mapsto 2^{-1}\,\phi_1(2x)$. In general, if

$$\phi_i(x) \xrightarrow{\ R\ } 2^{-\Delta_i}\phi_i(2x)\,,$$

then the number $\Delta_i$ is called the scaling dimension of $\phi_i$. For a lattice model, scale-invariance means scale-invariance of the mesh $\varepsilon \to 0$ limit, in which we rescale the fields $\phi_i$ by $\varepsilon^{\Delta_i}$. It is this limit that we actually observe on the macroscopic scale. See [7] for a superb exposition of scaling and renormalization.

Near a fixed point of $R$, we have much better chance of understanding what $R$ does. Many nearby theories will be attracted back to the fixed points by repeated applications of $R$. These should be put in the same universality class.

### A.3.3.

The critical Ising model should be scale-invariant. For $T \neq T_c$, there is a microscopic scale in the model set by the scale at which the signs exponentially decorrelate. At $T = T_c$ this becomes infinite and scale-invariance should appear. Currently, there is no mathematical

22 Like any fundamental idea in science, this one had many precursors in the work of many people. See e.g. [7, 22, 28, 31, 35–37] for various perspectives.

proof of this for $d = 3$ and it remains an important open problem. Numerical experiments [21] give

$$\Delta_{\sigma,3} = 0.518154\ldots$$

as the scaling dimension of the spin field in $d = 3$.

### A.3.4.

Importantly, for $d = 2$ and $d = 3$, this is not a Gaussian fixed point. A Gaussian random field is a generalization of a Gaussian process with $d$-dimensional time. For a Gaussian field, $\boldsymbol{\Phi} = \mathbb{R}^N$ and the functions $U_1$ and $U_2$ are quadratic. In suitable coordinates, the field thus becomes a superposition of many noninteracting Gaussian random variables. While certainly a very, very important part of probability theory and mathematical physics, Gaussian fields *do not* describe materials with nonlinear interactions.

By contrast, Theorem 3 implies that for $d \geq 4$, the critical Ising model is Gaussian and $\Delta_{\sigma,4} = 1$. It is a very difficult and important mathematical theorem to prove, and it undescore the crucial importance of dimensionality in statistical mechanics.

### A.3.5.

Which other microscopic models will fall into the critical Ising fixed point ? The crucial feature of the Ising model is the $\pm 1$ symmetry between the two possible values of $\sigma(\mathsf{v})$. One should expect that any potential $U_1$ which has two symmetric minima, like the function in (82)

$$U_1(\phi) = \qquad (82)$$

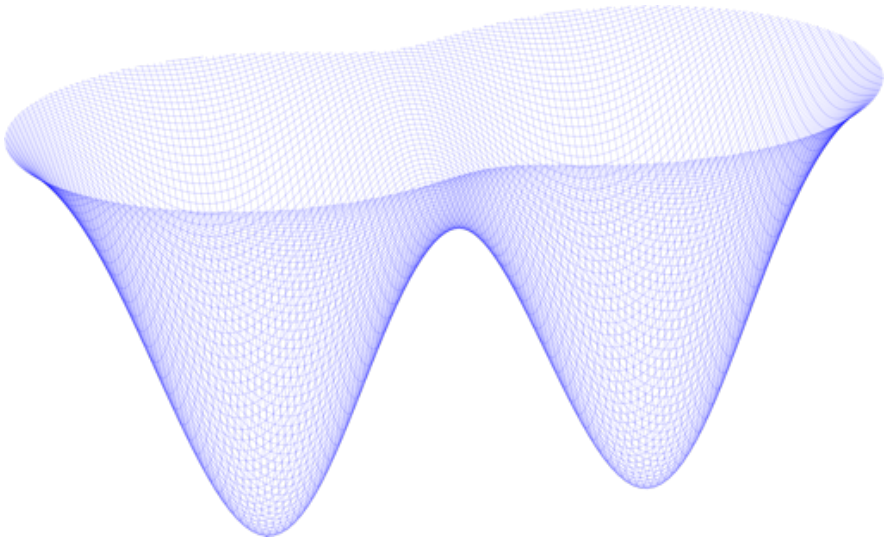

has the same critical scaling limit. For $T < T_c$, there should be Gibbs measures that prefer to stay close to one of the minima in (82), thus breaking the symmetry. For $T > T_c$ one expects a unique symmetric Gibbs measure.

While the details of these measures, and the values of the critical temperature $T_c$, will depend $U_1$ and $U_2$, the large-scale fluctuations of the $T = T_c$ measure (as captured by the the mesh $\varepsilon \to 0$ limit of the lattice model) should probably be universal and the same as for the critical Ising model. In other words, the critical Ising model should give the universal description of the transition between the broken and unbroken $\pm 1$ symmetry.

In light of Theorem 1, the reader may wish to contemplate what happens if the potential $U_1$ has 3 minima

$$U_1(\phi) = \tag{83}$$

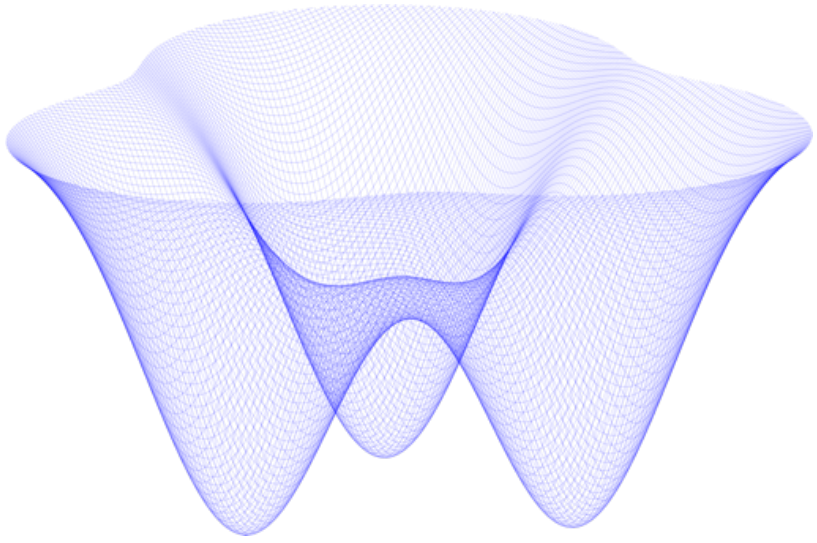

which can be permuted in all possible ways by the symmetries of the theory.

**Andrei Okounkov**

Andrei Okounkov, Department of Mathematics, University of California, Berkeley, 970 Evans Hall Berkeley, CA 94720–3840, okounkov@math.columbia.edu